\documentclass[11pt]{article}
\title{Clustering in Cell Cycle Dynamics with General Response/Signaling Feedback}

\author{Todd~Young\footnotemark[1] \footnotemark[4],  Bastien Fernandez\footnotemark[2], Richard~Buckalew\footnotemark[1], Gregory Moses\footnotemark[1], Erik Boczko\footnotemark[3]}

\usepackage{amsmath}
\usepackage{graphicx,color}
\usepackage{amsfonts}
\usepackage{amssymb}
\usepackage{verbatim}

\numberwithin{equation}{section}

\newtheorem{defn}{Definition}[section]

\newtheorem{cor}[defn]{Corollary}

\newtheorem{prop}[defn]{Proposition}

\newcommand{\D}[1]{{\mathbb#1}}% Doubled -Blackboard bold - caps only

\newcommand{\TT}{{\D{T}}}

\begin{document}

\maketitle

\centerline{To appear in J.\ Theoretical Biology {\bf  292} (2012), 103-115, doi:10.1016/j.jtbi.2011.10.002.}

{\small

\begin{abstract}
Motivated by experimental and theoretical work on autonomous oscillations in yeast, we analyze ordinary
differential equations models of large populations of cells with cell-cycle dependent feedback.
We assume a particular type of feedback that we call Responsive/Signaling (RS), 
but do not specify a functional form of the feedback. We study the dynamics and emergent behaviour of solutions,
particularly temporal clustering and stability of clustered solutions. We  establish the
existence of certain periodic clustered solutions as well as ``uniform''
solutions and add to the evidence that cell-cycle dependent feedback
robustly leads to cell-cycle clustering. We highlight the fundamental
differences in dynamics between systems with negative and positive feedback. For positive
feedback systems the most important mechanism seems to be the stability of individual isolated
clusters. On the other hand we find that
in negative feedback systems, clusters must interact with each other to reinforce coherence.
%For the form of feedback considered, there is a natural  constant $M$ that is the number
%of clusters that can exist without interaction. For $k \le M$, solutions consisting of
%$k$ noninteracting clusters in a positive feedback system are shown to be stable (neutrally), while
%for negative feedback such  solutions are unstable.
%Periodic solutions consisting of $k = M+1$ interacting clusters are shown to  be unstable
%for positive feedback and stable for negative. We study the dynamics of
%$k=2$ clusters completely.
We conclude from various details of the mathematical analysis that {\em negative} feedback is most
consistent with observations in yeast experiments.
\end{abstract}



}

\section{Introduction}

\footnotetext[1]{Mathematics, Ohio University, Athens, OH USA}
\footnotetext[2]{Centre de Physique Th\'{e}orique, CNRS - Aix-Marseille Universit\'{e}
Campus de Luminy, Marseille France}
\footnotetext[3]{Biomedical Informatics, Vanderbilt University Medical Center, Nashville, TN USA}
\footnotetext[4]{Corresponding author, {\tt youngt@ohio.edu}}

%
%\author{
%Todd R.~Young\footnote{Corresponding author, , Dept.~Math, Ohio Univ.}
%\and
%Bastien Fernandez\footnote{
%Centre de Physique Th\'{e}orique (UMR 6207 CNRS - Universit\'e Aix-Marseille~II - Universit\'e Aix-Marseille~I - Universit\'e Sud~Toulon-Var) CNRS Luminy Case 907, 13288 Marseille CEDEX 9, France}
%\and
%Richard Buckalew\footnote{Dept.~Math, Ohio University,
%Athens, OH USA}
%\and
%Gregory Moses\footnote{Dept.~Math, Ohio University,
%Athens, OH USA}
%\and 
%Erik M.~Boczko\footnote{Department of Biomedical Informatics,
%Vanderbilt University School of Medicine}
%%{\sf \small erik.m.boczko@vanderbilt.edu}
%}

\subsection{Background}

In this paper we consider simple dynamical models of the cell division cycle.
Specifically, consider a culture of $n$ cells, in which the progression
of the $i$-th cell is governed by the equation:
\begin{equation}\label{ode}
\frac{dx_i}{dt} = 1 + a(x_i,\bar{x}).
\end{equation}						
where $x_i$ is the position of the cell within the cycle and $\bar{x}$ denotes the state of all the
cells in the culture.  We will describe the dependence of $a$ on $x_i$ and $\bar{x}$ below.

Our primary motivation for this model is recent theoretical and experimental work 
on Yeast Autonomous Oscillations (YAO) (\cite{PNAS08}, \cite{tu}, \cite{chenz}, \cite{henson04}, \cite{klev3}), the periodic 
oscillations of physiologically relevant variables that have been reported and studied for over 
40 years~\cite{chenz,finn,keulers,kuenzi69,meyenburg69,muller,munch,murray,palkova,richard,satroutdinov,vonm,wang,xu}
(and many others). Different types of YAO have been called metabolic \cite{tu}, glycolytic \cite{bier} or respiratory \cite{henson04} oscillations. 
The control of oscillation and the regulation of yeast metabolism has been an important theme 
in the chemical engineering literature devoted to the efficient management of 
bioprocesses~\cite{buese,heinzle,hjortsO,fermentor,zhu}. These phenomena are of basic biological interest because they expose 
questions regarding the coordination of the cell cycle and metabolism, and interconnectedness of various cellular 
and genetic processes~\cite{pnas05,klev3}.
A correlation between YAO and the bud index was noted as early as \cite{kuenzi69,meyenburg69}. 
However, it seems that the link between YAO and the Cell Division Cycle (CDC) was obscured by the fact that the periods 
of YAO always shorter  than the CDC times (computed from dilution rate) and 
a relationship between YAO and the CDC seems to have been largely ignored.
However in \cite{klev3}, \cite{tu} and elsewhere, the correlation between YAO and CDC 
was again noted in genetic expression data.

\begin{figure}[th!]
\centering
\includegraphics[width=7cm,height=7cm]{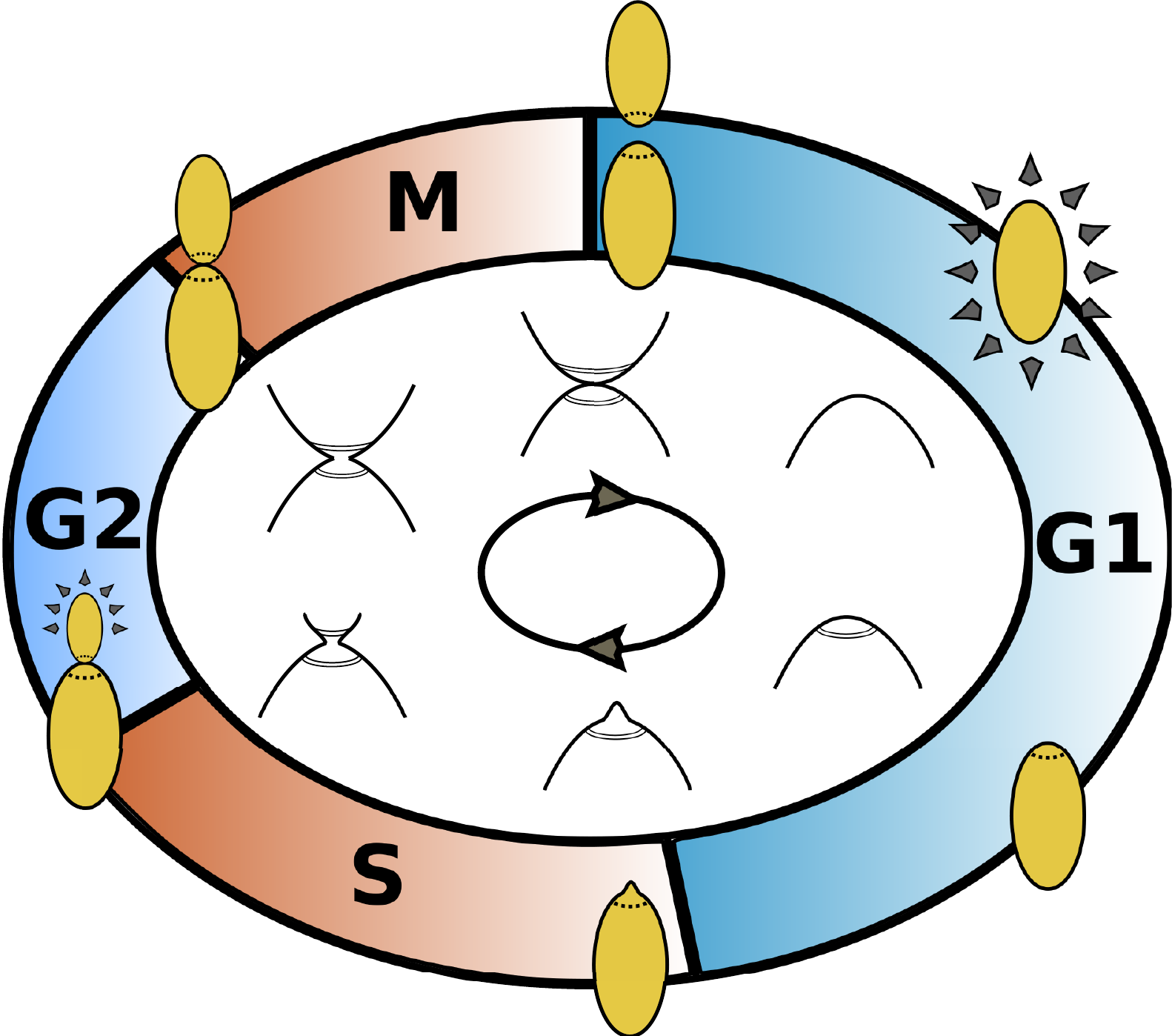}
\caption{Phases of the yeast cell cycle. The G1 phase begins following cell division.
The beginning of the DNA synthesis phase, S, coincides with budding. G2 is a second ``gap'' phase.
The M phase is characterized by narrowing  or ``necking'' between the parent and daughter cell;
it ends in cell division.  The hypothesized
$R$ region is the later portion of G1 and the signaling region $S$ is in the S phase.
That is, a large subpopulation of cells in the replicative S phase may promote or
inhibit progression of cells approaching the G1-S boundary. }
\label{cycle}
\end{figure}

In \cite{BSGY} and  \cite{PNAS08} the authors 
proposed {\em cell cycle clustering} as a possible explanation of the interaction between 
YAO and the CDC. Figure~\ref{cycle} roughly illustrates the arrangement of phases of the cell cycle of yeast.
We hypothesized that subtle feedback effects
on CDC progression could cause populations of cells to segregate into approximately
CDC-synchronized subpopulations. Experimental bud index data reported
in \cite{PNAS08} supported this picture.  In \cite{BSGY} the authors studied a few simple 
forms of (\ref{ode}) with the hypothesis that cells in one part of the CDC may influence  other cells in different parts of the CDC in different ways through various diffusible chemical products.
We hypothesized that a large subpopulation of cells in the critical S-phase might effect metabolism
production and the metabolites  may in turn
inhibit or promote cell growth in the later part of the G1 phase, thus setting up a feedback mechanism
in which YAO and CDC clustering are inextricably intertwined. 
We showed analytically and numerically that differential CDC feedback such as this  can robustly
cause CDC clustering in the models. By clustering we do not mean spatial clustering (cultures
that exhibit YAO occur in well-mixed bioreactors), but
groups of cells that are traversing the CDC in near synchrony.

\begin{figure}[!ht]
\centering
\includegraphics[width=15cm,height=10cm]{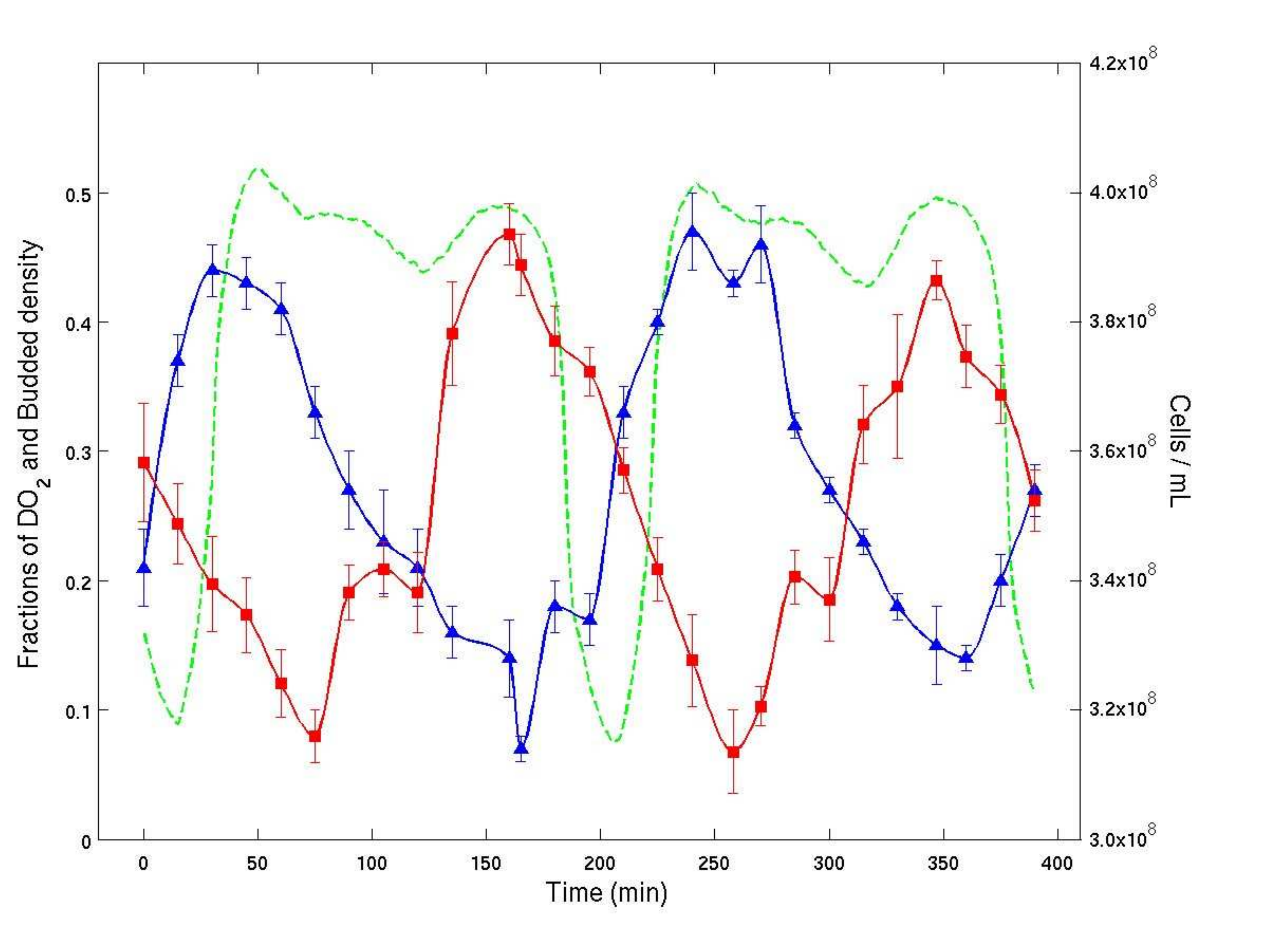}
\caption{Experimental time series from a continuous culture of budding yeast.
Dissolved $O_2$ percentage (green), bud index percentage (blue)
and cell density (red) are plotted versus elapsed time. The average cell
cycle period as calculated from dilution rate was about 400 minutes. The plot shows
clearly that the bud index (percentage of cells with buds from microscopy) and cell density
(by flow cytometry) 
are both synchronized with the oscillation in the level of dissolved $O_2$.
}
\label{ThePlot}
\end{figure}

Guided by these
mathematical results, we  verified the existence of clusters in two types of
oscillating yeast \cite{Stowers_comm} using both bud index and cell density data. Some of the 
measurements from those experiments are shown in Figure~\ref{ThePlot}.  First we note that the 
cell cycle period, as calculated by the dilution rate, is approximately 400 minutes, and
two $O_2$ oscillations occur during this period, suggesting that there may be two clusters.
Next, analyzing the figure, 
we see that approximately half of the cells are budding at times \( t = 50 \) and \( t = 250 \), 
while at \(t=170\) less than 10\% of the cells are budded.  Each budding event is accompanied by 
a decrease in density (no cells are dividing) and followed by a sharp increase in cell density as 
these budded cells proceed through division.  Note that at \(t=170\) since less than 10\% of the cells
are budded,  most of the cells must be in the G1 phase of the cycle. When the bud index
hits its next maximum at \(t=250\), approximately half of the cells must have budded. The other half of the cells
must at that time still be in the G1 phase. As the cells that are budded then divide and the cell density increases,
the cells that remained in G1 must still be in G1 since the bud index is again low. 
When the first group of cells has divided, the second group
has been in G1 for at least 200 minutes. The next rise in bud index then must be due to these cells, since
they have had time to mature, while the first group of recently divided cells clearly has not had time 
to reach budding again.  Thus, these experiments show conclusively the existence of two clusters and
that CDC clustering coexists with YAO.

\subsection{Modeling of the cell cycle and feedback}

In standard modeling the cell volume $v_i(t)$ is a proxy
for position in the cell cycle. This has justification for yeast in that
milestones in the cell cycle, such as the onset of budding, are closely associated
with volume milestones and thought to be causally related. Measurements show
that the growth of a single cell is roughly exponential, so a first order
approximation is that the volume of the $i$-th cell satisfies a linear differential equation
\begin{equation}\label{eqn:exp_growth}
\frac{dv_i}{dt} = c v_i.
\end{equation}
A frequent assumption on the growth rate $c$ is that it does not depend on $v_i$
{\em i.e.}\ it is independent of the cell's current state within the cycle and on other cells; it depends instead 
on the nutrients available and other environmental factors. 
Applying a logarithmic change of variables the growth law becomes $dx_i/dt = c$, and
by further normalizing both the coordinate $x$ and time, the cell cycle can be represented
by the unit interval $[0,1]$, and the equation of motion becomes $dx_i/dt = 1$.
In this simple model each cell reaches division (cytokinesis) at $1$ when it returns
(perhaps with its descendant cell) to $0$ and begins the cycle again.  Note that a
change to normalized coordinates does not depend essentially on the form of 
(\ref{eqn:exp_growth}); if there is no interaction between cells and cells never stop growing,
then one can change variables to the form $dx_i/dt = 1$, with $x_i(t) \in [0,1]$.

A much more general model (again using normalized coordinates) is (\ref{ode}).
We proposed to consider forms of (\ref{ode}) where the cells in one portion of the cell cycle,
$S$  for {\em signaling}, may influence the growth rate of cells in a preceding portion
that we term $R$ for {\em responsive}. For example the $R$ region may reside 
in the  later portion of the G1 phase and the signaling region $S$ may be the
biological S phase (see Figure~\ref{cycle}). This is philosophically justified
by the fact that the S phase is the most critical part 
of the CDC and the link between YAO and CDC may function to protect
the integrity of transcription \cite{chenz}. It is also known that growing yeast store carbohydrates, then 
liquidate them in the late G1 phase \cite{futcher}. The actual positions of the signaling
and responsive regions within the biological cell cycle play no role in the mathematical analysis. In the rest
of the paper we will use $S$ to denote the signaling region and on the few occasions when we
refer to the yeast's S phase we will do so explicitly.

Mathematically, the interval $[0,1]$, with the endpoints identified
is a circle. On this circle, we can specify a positive direction as associated with the increasing direction on $[0,1]$.
Distance between points $x$ and $y$ on the circle using these coordinates is given by the minimum of $|x-y|$ and $1-|x-y|$.
\begin{defn}
Consider $n$ cells whose coordinates are given by $x_i \in [0,1]$, $1$ identified
with $0$, and governed by an equation of the form (\ref{ode}). When a cell reaches $1$ it continues at $0$.
We call such a system a {\em RS-feedback system} if:
\begin{enumerate}
\item[(H1)] $R$ is an interval that directly precedes another interval $S$, 
                           i.e.\ the last endpoint of $R$ is the first endpoint of $S$,
\item[(H2)] $a(x_i,\bar{x})$ vanishes except when $x_i \in R$
      and there are some $x_j$ in $S$,
\item[(H3)] $0 < v_{min} \le 1+ a(x_i,\bar{x}) \le v_{max}$ for all $x_i$ and $\bar{x}$,
\item[(H4)] feedback is monotone, thus adding a cell to $S$ will increase
          the value $|a(x_i,\bar{x})|$ for $x_i \in R$, and,
 \item[(H5)] $a(x_i, \bar{x})$ is a smooth function for  $x_i $ in the interior or $R$ and each $x_j$ in the interior of $S$, $j \neq i$
                      and the one sided derivatives exist at the boundaries of $R$ and $S$.
\end{enumerate}
 By {\em positive feedback}  we mean $a$ is positive for $x_i \in R$ if there are one or more
$x_j$ in $S$.  We define {\em negative feedback} in the analogous way.
\label{defn:RS}
\end{defn}
For the sake of definiteness we will specify throughout the rest of the paper:
$$
     S = [0,s)  \quad \textrm{ and }   \quad  R = [r, 1) , \quad   0< s < r < 1.
$$
The final endpoint of $R$ is $1$, which corresponds to $0$,  the initial endpoint of $S$.
See Figure~\ref{model} below.

Note that our restriction that \emph{R} precedes \emph{S},  while  motivated by biological considerations,
is not the only possibility.  It is worth noting that in the reverse case, when \emph{R} follows \emph{S}, 
many of our results hold with the roles of positive and negative feedback reversed.  In particular, this 
is true for Propositions~\ref{pos_cluster_isolation},~\ref{pos_cluster_stability},~\ref{pos_two_clusters},~\ref{neg_cluster_stability},~\ref{neg_two_clusters}, and Corollary~\ref{neg_periodic_convergence}.

In~\cite{BSGY} we considered two idealized forms of (\ref{ode}) with threshold-triggered feedback and
found that with {\em either positive or negative linear feedback}, we robustly
observe clustering.  We also observed that the number of clusters formed is less
dependent on the form of the feedback than on the sizes of the signaling and responsive
regions, $S$ and $R$. One aim of the present paper is to
demonstrate that clustering typically occurs in a broad class of models (\ref{ode}), 
beyond the idealized situations investigated in \cite{BSGY}. The other aim is to study
key differences between positive and negative feedback in the models in order to
understand if one or the other is more favorable for clustering and thus more
likely to be the process behind experimental observations.

Note in Definition~\ref{defn:RS} that the number of cells in the culture
is fixed at  $n$. Thus we will consider neither death, harvesting or proliferation.
In the oscillation experiments we are modeling there are in fact proliferation and harvesting (and an insignificant
rate of death), but they approximately balance  when averaged over a cell cycle. Thus the expected number of
cells descended from a single cell at any future cycle is approximately $1$.  Further, in the model
we are considering, there is no distinction between the two cells resulting from a division and  to keep track
of both trajectories  would be redundant.

The differential equation in the general model (\ref{ode}) with RS feedback may have discontinuities.
Thus we need to consider uniqueness and global existence of solutions. First consider that 
the equations may be discontinuous only when a variable is at the boundaries of $R$ and $S$
{\em i.e.}\ at the hyper-surfaces given by $x_j = 0, s, r$ for some $j$. 
We obtain uniqueness since by (H3) any solution crosses a surface of discontinuity non-tangentially
with non-zero speed. Global existence follows from boundedness of the vector field (H3).

 A version of the model that we will study in sections 6 and 7  
is:
\begin{equation}\label{eqn:fmodel}
    \frac{dx_i}{dt}  = \begin{cases}
                      1, \quad \textrm{if} \quad x_i \notin R \\
                      1+  f(I)  , \quad \textrm{if} \quad x_i \in R
                     \end{cases}
\end{equation}
where 
\begin{equation}
I(\bar{x}) \equiv \frac{\#\{j:x_j \in S\}}{n}  = \frac{\#\{j:x_j \in [0,s) \}}{n}  ,
\label{eqn:I}
\end{equation}
\emph{i.e.}\ $I$ is the fraction of cells in the signaling region. 
The ``response function" $f(I)$ in (\ref{eqn:fmodel})  must satisfy 
$f(0) = 0$ and be monotone, but can be non-linear, for instance
sigmoidal (S-shaped). 
We will see  that models of this form, while fairly general, can be studied in some detail.
In \cite{FY10} we studied some cases where $f$ is a linear function.

Understanding the  CDC at the genetic and biochemical level is a topic of intense interest
and  progress has been made in identifying the agents and  the nature
of relationships between them \cite{boye,collier,futcher,singer,spellman,tvegard}.
Our approach uses a ``caricature" of the cell cycle, rather 
than detailed modeling, and this simplification demands justification.  First, we wish
to deal with individual cells in a population-wide phenomena. If details within each cell are considered,
then the dimensions of the phase space would be extremely large and results would be 
difficult to obtain. Second,  our understanding 
of the details of the cell cycle and its relations
with other processes is not complete and even if the general nature of
relationships were well-understood, the resulting set of differential equations would
contain many parameters, {\em e.g.}\ rate constants, that could only be estimated. 
With our simplified model which is based on biological insight, 
we hope to obtain general principles that will inform further detailed investigations. 

The approach in part of this work is basically that of  ``phase oscillator" models, e.g.~Kuramoto equations, in which
details of each individual actor are projected onto a simple phase space and emergent population behaviours are studied.
In fact if $f(I)$ is linear then our model can be put into the Kuramoto form by integrating over
the cell cycle for each pair of cells and adding the effects (see \cite{kura} p. 65-67). 
This derivation fails in general for RS models (\ref{ode}) or  (\ref{eqn:fmodel}) 
since the effects of cells in general are not additive.

We note here another modeling simplification that we have implicitly made; namely one might more accurately
model the feedback term as $a(x_i,z)$ where $z$ is a vector variable representing all substrate
factors that contribute to growth rate and $z$ itself is coupled with $\bar{x}$ \cite{hjo}.
Dropping the $z$ variable can be justified if the time-scale of the dynamics of this variable
is significantly shorter than the time-scale of the CDC.

%%%%%%%%%%%%%%%%%%%%%%%%%%%%%%%%%%%%%%%%%%%%%%%%%%%%%%%%%%%%%%%%%%%
%%%%%%%%%%%%%%%%%%%%%%%%%%%%%%%%%%%%%%%%%%%%%%%%%%%%%%%%%%%%%%%%%%%
%%%%%%%%%%%%%%%%%%%%%%%%%%%%%%%%%%%%%%%%%%%%%%%%%%%%%%%%%%%%%%%%%%%
%%%%%%%%%%%%%%%%%%%%%%%%%%%%%%%%%%%%%%%%%%%%%%%%%%%%%%%%%%%%%%%%%%%

\section{Clusters, Gaps and Isolation}

In this section we begin to study the existence and stability of
periodic  ``clustered'' solutions for both positive and negative
$RS$-feedback systems (\ref{ode}) and in later sections we point out crucial differences between these
two types of feedback. Reducing to the study of clustered solutions is 
of practical interest since clusters appear in experiments with YAO. It also limits
the dimensions of the problem to a manageable size. This strategy has proven indispensable 
in many fields; for instance in 
fluid dynamics, insight is obtained by studying finite dimensional
vortex equations rather than the full Navier-Stokes partial differential 
equations \cite{newton01_vortex}.

\begin{defn}
By a {\em cluster}, we will formally mean a group of cells that are completely synchronized
in the CDC. 
\end{defn}

Note that RS-feedback systems as defined have the symmetry of globally coupled networks with
identical nodes; namely, the vector field is equivariant with respect to the group of
permutations of coordinates. This symmetry implies that any cells that initially
share the same phase keep the same phase as time evolves.
The simplest trajectory consists in taking all cells initially in the same phase. We have
a single cluster $C$ (synchrony) that generates a periodic solution
that runs at velocity 1 around the circle ({\em i.e.}\ $x_i(t)=x_i(0)+t$ $\textrm{mod } 1$ for
every $i$ and all $t>0$).

\begin{defn}
 By the {\em gap} between two clusters or cells at $x_{i-1}$ and $x_i$ we mean  the open 
 interval on the circle from $x_{i-1}$ to $x_i$, in the direction of the flow that contains no
 cells and has width $w_i = x_i - x_{i-1}$. (This can be made precise using a lift  to the real line.)
\end{defn}
Note that if there are only two clusters $x_0$ and $x_1$ in the system, then there are two gaps between them
and each gap has a width.

\begin{figure}[!t]
\centering
\includegraphics[height=8cm,width=11cm]{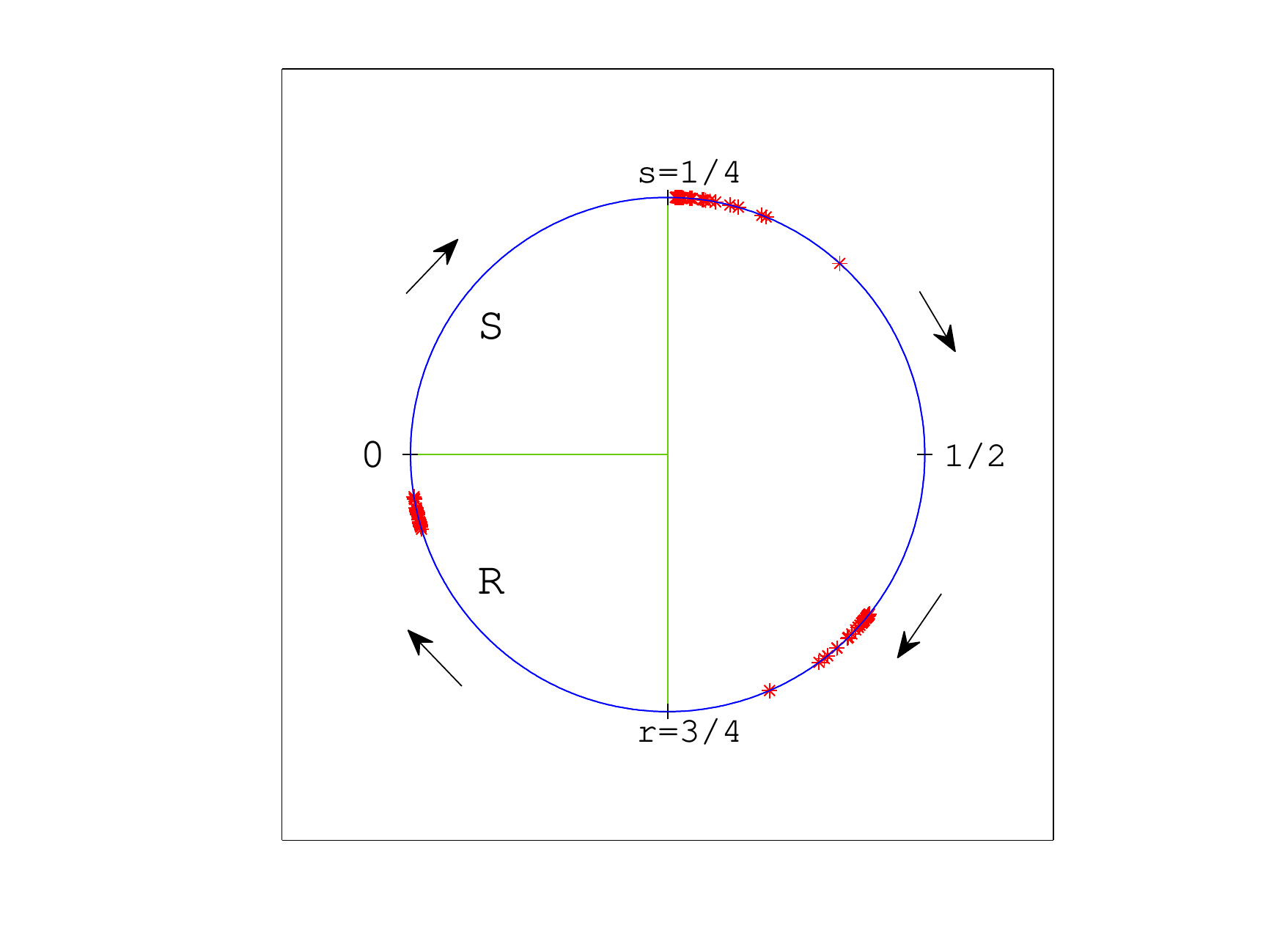}
\caption{Our coordinate representation of the cell cycle and (weakly clustered) groups of cells from a negative feedback
simulation with $n = 200$ and parameter values $s=.25$ and $r = .75$. Positions of individual cells
are denoted by red asterisks.
In this coordinate system  the $S$ region is the interval $[0,.25)$ and the
$R$ region is $[.75,1)$, where $1$ is identified with $0$.
}
\label{model}
\end{figure}

It follows from our assumed coordinates that if two clusters are in $S$ and $R$ then the 
distance between them on the circle is less than $|R|+|S| = 1-r+s$. Here $|R|$ denotes
the length of the interval $R=[r,1)$, which is $1-r$ and $|S|=s$ denotes the length of the
interval $S=[0,s)$. See Figure~\ref{model}.
We say that a cluster of cells is {\em isolated} if there are gaps between
the cluster and  any other cells on either side of length at least $|R|+|S|$ and 
{\em strictly isolated} if the widths of gaps are more than $|R|+|S|$. This terminology is motivated by the fact that
strictly isolated clusters cannot exert feedback on cells outside the cluster, or have feedback exerted upon them from
outside. 
While we consider only clustered solutions in the strictest sense, in real
cultures individual cell differences  will lead to a weaker form of clustering.
 For clarity we will refer to such a weakly clustered subset of cells as a {\em group} of cells.
In Figure~\ref{model}  the solution in the simulation has formed groups.

The following definition will play a large role in the analysis of the model.

\begin{figure}[h!]
\centering
\includegraphics[width=11cm,height=9cm]{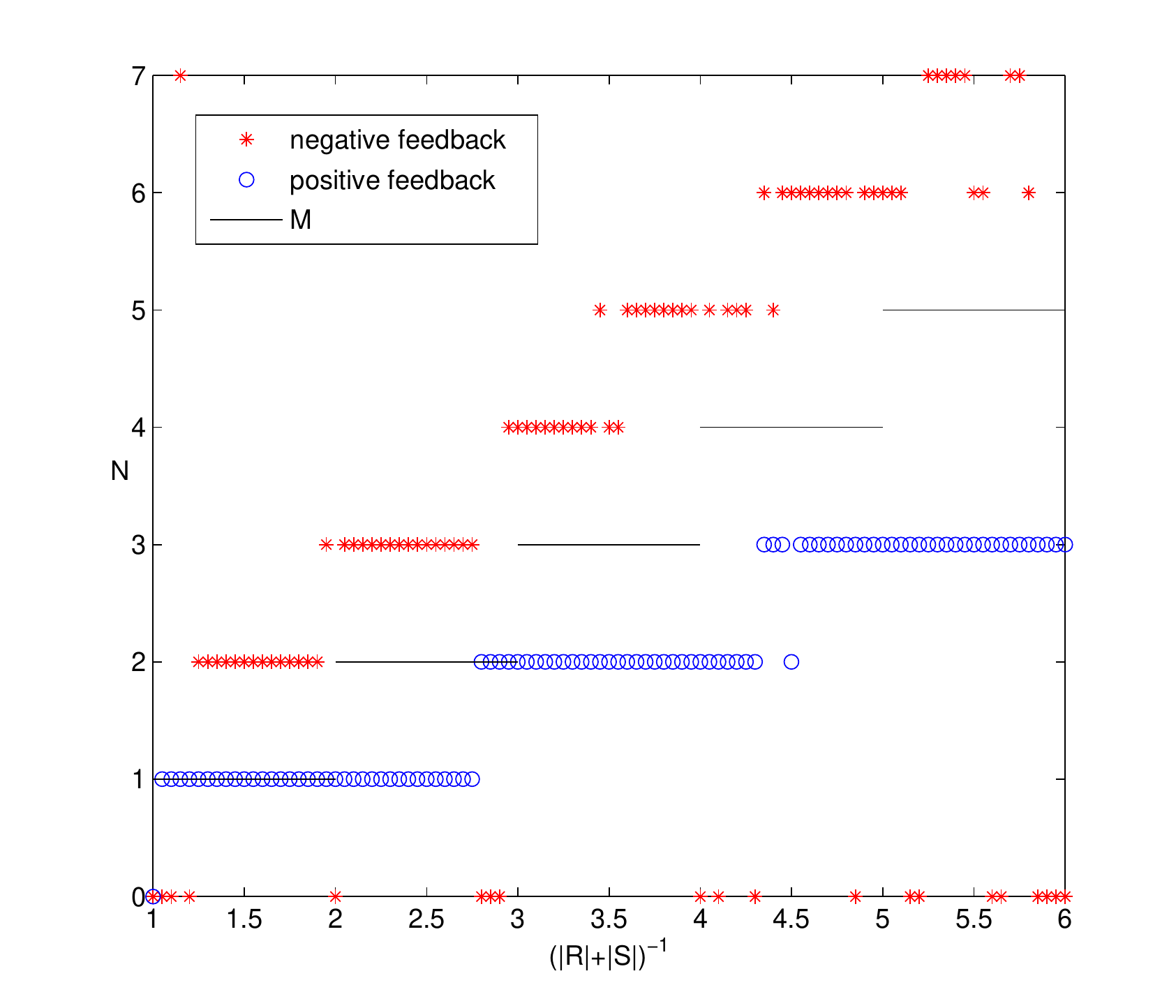}
\caption{The number of clusters that form in simulations compared with $M = \lfloor (|R|+|S|)^{-1} \rfloor$, the
               maximum number of isolated clusters given $R$ and $S$.
}
\label{number_clusters}
\end{figure}

\begin{defn}
Define 
\begin{equation}
M \equiv \lfloor (|R|+|S|)^{-1} \rfloor.
\label{eqn:M}
\end{equation}
 M is the maximum number
of  isolated clusters that can simultaneously exist, given the sizes $|R|$ and $|S|$ of $R$ and $S$. 
\label{def:M}
\end{defn}
Here $\lfloor x \rfloor$ denotes the floor function, that is, the greatest integer less than or equal to $x$ (e.g.\ $\lfloor 2.1 \rfloor = 2$). 
In Figure~\ref{number_clusters} when $(R|+|S|)^{-1}$ is in the range $[3,4)$, $M$ is $3$.

\begin{prop} For any RS feedback model  and any positive integer $k \le M$ there exist periodic solutions 
consisting of $k$ isolated clusters
that do not interact.
\end{prop}
{\sc Proof:}
For $k\leq{M}$ consider the solution with initial conditions:
$x_{0}=0$, $x_{1} = \frac{1}{k}$, ..., $x_{k} = \frac{k-1}{k}$. 
We claim that this is such a periodic solution.
Notice from the definition of $M$ that  
$$
M \le (|R|+|S|)^{-1} 
$$
and so
$$
 \frac{1}{k}  \ge \frac{1}{M} \ge |R|+|S|  = 1-r+s  .
$$
Since the distance between any two consecutive clusters  $x$ and $y$
is initially $d(x,y) = 1/k$, no two clusters can be in $R \cup S$ simultaneously. 
Thus no feedback will occur and thus the distance between clusters will not
change. 

By the same reasoning  any initial condition of $k \le M$ clusters, where all
pairs of clusters satisfy $d(x,y) \ge |R|+|S|$, will also lead to a periodic solution where all
clusters move indefinitely with speed 1.
\hfill $\Box$

Conversely, if more than M clusters exist, then at least two of them are within a distance $|R|+|S|$ of each other,
and while the first of these clusters lies in the signaling region, it will exert feedback on the second cluster for
a non-empty interval of time.

Note that a solution consisting of strictly isolated clusters  can be at most
neutrally stable (not asymptotically stable) since moving a cluster to the left or 
right still produces an isolated cluster.

In Figure~\ref{number_clusters} we plot the results of numerical simulations
which compare the number of clusters that  formed with the maximum number
of possible isolated clusters. The number of cells $n$ was 5000
and the model was taken to be that in equation (\ref{eqn:fmodel}) with $f$ linear.  Specifically, for
the positive feedback simulations $f(I) = .6 I $, and $f(I) = -.6 I$ in the negative feedback simulations.
There was also stochastic noise added to the equations with level $\sigma = 10^{-6}$ in order
to demonstrate robustness.
In the plot the $x$-axis represents $(|R|+|S|)^{-1}$ which was varied from $1$ to $6$
in one hundred increments. In the simulations $|R|$ and $|S|$ were taken to be equal.
The system was integrated for up to 100 cell cycles (with 50 steps per cycle) to check for the formation of
clusters (but clusters usually formed long before the 100th cycle). To test for clustering, we
produced histograms with 120 bins and visually inspected these for clustering.
If clusters did not clearly form for a given value of $(|R|+|S|)^{-1}$, then that data point
is plotted on the graph as a zero value for $N$.

The most striking feature of the plots in Figure~\ref{number_clusters} is that
for positive feedback the number of clusters formed is always less than or equal
to $M$, but for negative feedback the number of clusters formed is always
greater than $M$. Notice also that positive feedback always produced clustering,
but this was not the case for negative feedback. For negative
feedback there are parameter values where
no clusters form. Finally, it is worthy of note that for negative feedback, there are no
occurrences of one-cluster solutions. Analysis in the next sections will
shed  light on these observations.

Besides solutions consisting of $k\leq M$ isolated clusters,
RS-feedback systems have other periodic solutions.
One of these consists of all $n$ cells spread along the cycle as uniformly as possible.
We will define a {\bf uniform solution} to be a
trajectory for which the coordinates satisfy the following relation for some time $d>0$:
\begin{equation}\label{uniform}
    x_i(d)=x_{i+1}(0) \quad \textrm{for all} \quad  i=0, \ldots, n-2,
           \quad \textrm{and} \quad x_{n-1}(d)=x_0(0)  \mod 1.
\end{equation}
%In particular, due to the cyclic permutation symmetry of the system, the configuration of a uniform
%solution at time $t+d$ ($t>0$) coincides, up to a relabeling of indices, with the
%configuration at time $t$. 
Since the velocity of cells in the complement of $R$ is precisely $1$,
it follows that for such a solution the cells in $R^c$ will be uniformly
distributed with inter-cell distance $d$. In the process of the construction of uniform solutions
 we also produce many other periodic $k$ cluster solutions. 
 
 Suppose that $k$ divides $n$ and
 and $k$ sets of $n/k$ cells are initially synchronized. Then we may greatly reduce the dimensions
 of the differential equations by considering only the positions of the $k$ clusters which we may denote
 by $\{x_0(t), x_1(t), \ldots , x_{k-1}(t) \}$. In particular we will prove that there always exists a solution
 of $k$ clusters that satisfy:
 \begin{equation}\label{cyclick}
     x_i(d)=x_{i+1}(0) \quad \textrm{for all} \quad i=0, \ldots, k-2,
                 \quad \textrm{and} \quad x_{k-1}(d)=x_0(0) \mod 1.
\end{equation}
We will refer to such solutions as {\bf cyclic $k$ cluster solutions}.

\begin{prop}
 There exists a uniform solution of any RS-feedback system. If $k$ is a divisor of $n$, then
 a cyclic $k$ cluster solution exists consisting of $n/k$ cells in  each cluster.
\end{prop}
We defer the proof of this proposition until Section 5.2.  We observe that we do not have uniqueness of the uniform and cyclic clustered solutions.

For $n$ large, as in the application in mind, then we expect $k$ cluster solutions for $k <<n$ to exist even if
$k$ does not exactly divide $n$. For example hyperbolicity (linear stability or instability) of the
$k$ cyclic solutions when
$k$ divides $n$ would imply that such solutions exist for all $n' \approx n$.

The uniform solution has an analogue in PDE models of the cell cycle which we note in Appendix~A.

%%%%%%%%%%%%%%%%%%%%%%%%%%%%%%%%%%%%%%%%%%%%%%%%%%%%%%%%%%%%%%%%%%%%%%%%%%%%%%%%%%
%%%%%%%%%%%%%%%%%%%%%%%%%%%%%%%%%%%%%%%%%%%%%%%%%%%%%%%%%%%%%%%%%%%%%%%%%%%%%%%%%%
%%%%%%%%%%%%%%%%%%%%%%%%%%%%%%%%%%%%%%%%%%%%%%%%%%%%%%%%%%%%%%%%%%%%%%%%%%%%%%%%%%
%%%%%%%%%%%%%%%%%%%%%%%%%%%%%%%%%%%%%%%%%%%%%%%%%%%%%%%%%%%%%%%%%%%%%%%%%%%%%%%%%%

\section{Positive Feedback Systems}

In a RS model with positive feedback, first notice that a group of cells that is 
isolated will remain isolated. Further, positive feedback has a focusing effect on an isolated group.
\begin{prop}\label{pos_cluster_isolation}
  In a general RS model (\ref{ode}) with positive feedback, suppose that a solution has a
  gap between two adjacent cells $x_{i-1}$ and $x_i$ of width greater than or equal to $|R|+|S|$. 
  Then the width of this gap will never decrease. In particular, an isolated 
  group or isolated cluster will remain isolated indefinitely.
\end{prop}
If there are only two clusters in the system, then  this proposition applies to either of the two gaps that 
has width at least $|R|+|S|$.

{\sc Proof.} 
Suppose that two consecutive cells $x_{i-1}$ and $x_i$ are separated by a gap of width $w_i \ge |R|+|S|$.
Since feedback is positive, the cell $x_i$ always moves at a speed of at least $1$. The cell speed of
$x_{i-1}$, being governed by (\ref{ode}), will be exactly $1$ whenever $w_i$ is greater than or equal to $|R|+|S|$,
since there will be no cells in $S$ when $x_{i-1}$ is in $R$.
We then find that the time derivative of $w_i$ is non-negative when $w_i \ge |R|+|S|$. 
Therefore, the gap can never decrease. It follows immediately that a group or cluster of cells that is
isolated will remain isolated. \hfill $\Box$

\begin{prop}\label{pos_cluster_stability}
  In a general RS model (\ref{ode}) with positive feedback, suppose that a group of cells
  with width $w$ less than $|R|+|S|$ is isolated. Then  the width of the group 
  will converge to zero as $t \rightarrow +\infty$.
  \end{prop}
{\sc Proof.} 
By the previous proposition the group will remain isolated for all future time and we may consider this group of
cells as a decoupled sub-system. Without loss of generality we may renumber the cells in the group
so that they have coordinates: $x_1, x_2, \ldots, x_\ell$, (ordered in the direction of the flow).
 For $i=1, \ldots, \ell -1$, denote by $w_i = x_\ell - x_i$ the width of the interval from $x_i$ to $x_\ell$.
By assumption $w_{\ell-1} \le w_i \le w_1 < |R|+|S|$.
Observe that this condition ensures that each $x_i$, $i = 1, \ldots, \ell -1$ 
will experience some acceleration every time it passes through $R$, unless $x_i = x_\ell$. This implies
that if $w_i(t)$ is non-zero, then $w_i(t)$ will decrease each time the group passes through $R$ and $S$.

Since the group will remain isolated, the cell $x_\ell$ will always move with speed $1$ and $w_i(t)$ will never
increase. Since each $w_i(t)$ is non-increasing and bounded below by $0$, it must have a limit $w_i^\infty$.
 Now consider a solution with an initial condition such that $x_i(0) = x_\ell(0) - w_i^\infty$ for $1 \le i < \ell$. 
 By a standard argument using continuous dependence of solutions on initial conditions, 
 for this solution each $w_i(t)$ will be identically $w_i^\infty$. 
 This  implies that $w_i^\infty$ must be zero, since, from above,
 a non-zero $w_i(t) < |R|+|S|$ must decrease  during each cycle.
\hfill $\Box$

In the next result we discuss stability, for which we need the concept of distance and neighborhoods
in phase space, which for the  models we are considering is the $n$-torus, $\TT^n$, where
$n$ is the number of cells.   Note that on $\TT^n$ there is a natural metric (distance) defined by the
maximum of the (mod 1) coordinate differences.

\begin{prop}
In a RS model (\ref{ode}) with positive feedback, the set of strictly isolated cluster solutions is
locally asymptotically stable. A solution consisting of $k \ge 2$ strictly isolated  clusters
 is neutrally stable (stable, but not asymptotically stable) within the set of solutions with $k$ clusters.
\end{prop}

{\sc Proof.}
First observe that an $\epsilon$-neighborhood of a configuration consisting of isolated  clusters
consists of groups of cells within $\epsilon$ of the original clusters. If the original clusters
are strictly isolated, then we may make $\epsilon$ small enough that the groups
are also strictly isolated. By  
Propositions \ref{pos_cluster_isolation}~and~\ref{pos_cluster_stability}
each of these groups will remain isolated and converge to a  cluster. Thus a solution starting at any initial condition
within a neighborhood of the set of strictly isolated clusters will asymptotically approach the set.

The second part of the claim follows since, if a strictly isolated cluster is moved a small distance,
then it is still strictly isolated. Thus a small perturbation of a solution consisting of $k$ strictly isolated  
clusters also will be a solution consisting of $k$ strictly isolated clusters. The distance between the
two solutions will remain constant for all future time and thus they are stable, but not asymptotically stable.
\hfill $\Box$

Points near the set of isolated cluster solutions will converge to the set, but  individual solutions
are only neutrally stable with respect to perturbations inside the set.

Since the set of clustered solutions is locally stable, it must have a basin of
attraction and it is natural to ask how big the basin is. In simulations for
positive feedback systems, the basin seems to include almost all initial conditions.
In the next proposition we see that the basin of attraction extends far beyond
a small neighborhood of the set.

%\begin{prop}
%In a RS model with positive feedback, suppose that there are two clusters
%that are more than $|R|+|S|$ apart from each other and the leading cluster
%is isolated ahead. Suppose that there is another cluster in the gap between
%these clusters, within $|R|+|S|$ of the cluster ahead of it, and that the number of cells in the leading cluster is
%sufficiently large and the cluster in the middle consists of sufficiently few cells.
%Then the middle cluster will asymptotically converge to the leading cluster.
%\end{prop}
%\noindent

\begin{prop}\label{pos_two_clusters}
 Suppose that a solution $\bar{x}(t)$ in a RS model (\ref{ode}) with positive feedback has at least one
 gap of width greater than or equal to $|R|+|S|$. Then the solution will converge to a periodic solution
 consisting entirely of isolated  clusters.
\end{prop}
\noindent
{\sc Proof.}
We will call a gap {\em large} if its width is greater than or equal to $|R|+|S|$.
At time $t_0$, the cells may be grouped into a minimum number of groups in which
there are no large internal gaps.  The number of
such groups is the same as the number of large gaps. Note that
each such group is isolated, and the number of such groups cannot be larger
than $M = \lfloor (|R|+|S|)^{-1} \rfloor$. Consider one such group. Since it contains
no large internal gaps and it is isolated, during passage through $R$ the last cell must
be accelerated by the presence of at least one cell in $S$ and so its speed is
sometimes greater than $1$. On the other hand, since the group is isolated it
will remain isolated by Proposition 3.1, and the lead cell will travel indefinitely with 
speed $1$ by the same argument that appears in the proof of that proposition.

Thus during one passage through $R$ the distance between the lead cell and the final cell
in the group must decrease. If this group continues to have no large internal gaps, 
then it follows that the width of the group will continue to decrease. By
an elementary argument (as for Proposition 3.2), the width will converge to zero; in other words
the group will converge to a  cluster. Otherwise, if a large internal gap
develops then the cells that are separated by the gap will be isolated from each
other and thus form two isolated groups. When this occurs the number of isolated
groups will increase. Since the number of isolated groups is bounded above,
large internal gaps may form only a finite number of times and thus eventually we have
a fixed number of groups that never develop large internal gaps and each of these converges
to an isolated  cluster.
\hfill $\Box$

%%%%%%%%%%%%%%%%%%%%%%%%%%%%%%%%%%%%%%%%%%%%%%%%%%
%%%%%%%%%%%%%%%%%%%%%%%%%%%%%%%%%%%%%%%%%%%%%%%%%%
%%%%%%%%%%%%%%%%%%%%%%%%%%%%%%%%%%%%%%%%%%%%%%%%%%
%%%%%%%%%%%%%%%%%%%%%%%%%%%%%%%%%%%%%%%%%%%%%%%%%%

\section{Negative Feedback Systems}

The key observation is that for negative feedback, isolated clusters are not
stable. This is because as a group of cells crosses the R-S boundary
all cells  of the group are delayed except the lead cell, which moves with
unit velocity, causing the group to spread.

\begin{prop}\label{neg_cluster_stability}
 In a RS model (\ref{ode}) with negative feedback, a solution consisting of strictly isolated  clusters is locally unstable.
\end{prop}
\noindent
{\sc  Proof.}
Denote by $x^*(t)$ a solution consisting of strictly isolated  clusters.
First observe that under the condition of strict isolation the gaps between clusters are all
larger than $|R|+|S|$, and so  any sufficiently small perturbation of the  clusters
consists of groups that are still isolated.
In any neighborhood of any configuration with an isolated  cluster, there is
a configuration where $x_i \neq x^*_i(0)$ for any $i$. 
Note in fact that a local coordinate system in a neighborhood of
$x^*(0)$ is $\{x_i - x^*_i\}$, $0 \le i \le n-1$.
Now let $x(t)$ be a solution with initial condition $x(0)$ that differs from $x^*(0)$
in only the $i$-th coordinate.  When this cluster passes
through the boundary from $R$ to $S$ the separation between $x_i$ and the rest of
the cluster to which it belongs will increase. If we further let the perturbation be sufficiently
small, then the cluster will remain isolated.

Now recall the definition of stability: given any $\epsilon>0$, there exists $\delta >0$ such
that any solution $x(t)$ starting within a $\delta$ neighborhood of $x^*(t)$ will remain indefinitely
within an $\epsilon$ neighborhood of $x^*(t)$.  Let $\epsilon_0$ be the largest $\epsilon$
so that any $x(0)$ within an $\epsilon$ neighborhood of $x^*(0)$ will consist of isolated groups and 
let $\epsilon = \epsilon_0/2$. If $x(0)$ is as in the previous paragraph and is arbitrarily
close to $x^*(0)$ then the distance between  $x(t)$ and $x^*(t)$ will continue to increase
on each unit time interval as long as $x(t)$ continues to consist of isolated groups.
Therefore it follows that $x(t)$ will eventually be outside of an $\epsilon_0/2$ neighborhood
of $x^*(t)$. Thus $x^*(t)$ is not stable.
\hfill $\Box$

Note that  we have not proved linear instability (when the derivative of the mapping at the periodic point
is greater than one in absolute value which implies that nearby orbit are repelled exponentially). 
Linear instability can be shown for systems of the form
(\ref{eqn:fmodel}).

It follows that in order for clusters to remain coherent under small perturbations
in negative feedback, they must not be isolated, {\em i.e.}\ the gaps between them
must be less than $|R|+|S|$, and so the number of stable clusters must be at least
$M+1$. This is clearly confirmed in simulations. Further it seems that sometimes $M+2$
is a stable number of clusters. 
We will show in \S 7 that a $k = M+1$ cluster cyclic solution is stable for 
negative feedback of the form (\ref{eqn:fmodel})
and for some open sets of parameter values.

In the following proposition we see that interacting clusters tend to spread out from each
other as far as possible.
\begin{prop}\label{neg_two_clusters}
  In a RS model (\ref{ode}) with negative feedback suppose that two clusters are within $|R|+|S|$ of each other,
but are isolated from other cells (non-empty).  If they remain isolated from  other cells, 
then the gap between the two clusters will increase and converge to $|R|+|S|$. In the case that the two clusters contain all the cells in
the  system, if one gap has width less than $|R|+|S|$ and the second gap has width $\ge |R|+|S|$, then
as long as the width of the second gap remains greater than $|R|+|S|$ the first gap will increase and approach $|R|+|S|$.
\end{prop}
\noindent
{\sc Proof.}
If the gap width (or smaller gap width in the case of only two clusters) is less 
than $|R|+|S|$ then each time the second cluster passes through $R$, the first 
cluster will be in $S$ for a non-empty interval of time. During this time interval, the cluster in R will 
experience deceleration, and the width of the gap will increase during the passage through $R$.
If the cluster pair remains isolated, then the distance will be preserved through the rest
of the cycle. Thus the distance between clusters will increase during each cycle. The distance
is bounded above by $|R|+|S|$ and so, by a standard argument, the sequence of 
distances thus generated will converge to $|R|+|S|$.
\hfill $\Box$

\begin{cor}\label{neg_periodic_convergence}
 Suppose that there are $k \le M$  clusters in a RS system (\ref{ode}) with negative feedback.
Then the solution will converge to a periodic  clustered solution with isolated clusters.
\end{cor}

%%%%%%%%%%%%%%%%%%%%%%%%%%%%%%%%%%%%%%%%%%%%%%%%%%%%%%%%%%%%%%%%%%%%%%%%%%%%%%%%%%
%%%%%%%%%%%%%%%%%%%%%%%%%%%%%%%%%%%%%%%%%%%%%%%%%%%%%%%%%%%%%%%%%%%%%%%%%%%%%%%%%%
%%%%%%%%%%%%%%%%%%%%%%%%%%%%%%%%%%%%%%%%%%%%%%%%%%%%%%%%%%%%%%%%%%%%%%%%%%%%%%%%%%
%%%%%%%%%%%%%%%%%%%%%%%%%%%%%%%%%%%%%%%%%%%%%%%%%%%%%%%%%%%%%%%%%%%%%%%%%%%%%%%%%%

\section{Dynamics of clusters  via return maps}

\subsection{The return map for a clustered system}

Let us continue to consider the most general model (\ref{ode}) with RS feedback. 
Let a population of $n$ cells be organized into $k$ equal  clusters and
let the clusters be labeled by a discrete index $i\in\{0,\dots,k-1\}$ so that
$x=\{x_i\}_{i=0}^{k-1}$ represent clusters of $n/k$ cells each. 
One can assume that all coordinates $x_i(0)$ of the $k$ clusters are initially well-ordered as
\[
0 = x_0(0) \leq \ldots \leq x_i(0) \leq\ldots \le 1,  \quad i=1, \dots, k-1.
\]
This ordering is preserved under the dynamics (this can be well-defined using the orientation
of the circle and $x_0(t)$ as a moving reference point). Moreover, the first coordinate $x_0$ must
eventually reach 1, {\em i.e.}\ there exists $t_\text{R}$ such
that $x_0(t_\text{R})=1$.
Thus the set $x_0 = 0$ defines a Poincar\'e section for the dynamics and the mapping
\[
(x_1(0),x_2(0),\dots,x_{k-1}(0)) \mapsto (x_1(t_\text{R}),x_2(t_\text{R}),\dots,x_{k-1}(t_\text{R}))
\]
defines the corresponding return map.

\begin{figure}[h!]
\centering
\includegraphics[width=11.5cm]{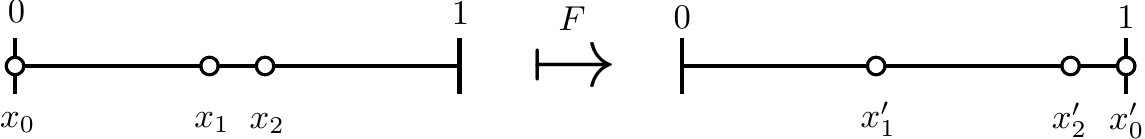}
\caption{An illustration of the map $F$ with $k=3$; $F(x_1, x_2) = (x'_1, x'_2)$.}
\label{F}
\end{figure}

Starting from $t=0$, compute the time $t_1$ that $x_{k-1}$ needs to reach 1 and compute 
the location of the remaining clusters at this time. Define  a map $F$ by
$$
  F : (x_1(0),x_2(0), \dots,x_{k-1}(0)) \mapsto (x_0(t_1),x_1(t_1),\dots,x_{k-2}(t_1)).
$$
Notice that $x_0(t_1)=t_1$ by assumption on $x_0(0)$.  An illustration of $F$ in the case \( k = 3 \) is given in Figure~\ref{F}.

Now the time $t_1+t_2$ that $x_{k-2}$ needs to reach 1, together with the population
configuration at $t=t_1+t_2$, follow by applying $F$ to the configuration
$(x_0(t_1),x_1(t_1),\dots,x_{k-2}(t_1))$. By repeating the argument,
the desired return time $t_\text{R}$ is given by $t_\text{R}=t_1+t_2+\dots +t_k$
and the desired return map is $F^k$. Therefore, to study the dynamics one only has to 
understand the first map $F$. 

Note the Kuramoto derivation (in which one averages over a cell cycle for each pair of cells then 
adds the effects (\cite{kura} p. 65-67)) is a specific way to calculate $F^k$ that is applicable 
when the feedback is additive.

We will first consider general properties of the map $F$ and then we will compute and analyze it
for  feedback of the form (\ref{eqn:fmodel}) in the simplest case, $k=2$. We emphasize that the case 
$k=2$ is also perhaps most important since it is the first to have been directly observed.

%%%%%%%%%%%%%%%%%%%%%%%%%%%%%%%%%%%%%%%%%%%%%%%%%%%%%%%%%%%%%%%%%%%%%%%%%%%%%%%%%%
%%%%%%%%%%%%%%%%%%%%%%%%%%%%%%%%%%%%%%%%%%%%%%%%%%%%%%%%%%%%%%%%%%%%%%%%%%%%%%%%%%
%%%%%%%%%%%%%%%%%%%%%%%%%%%%%%%%%%%%%%%%%%%%%%%%%%%%%%%%%%%%%%%%%%%%%%%%%%%%%%%%%%
%%%%%%%%%%%%%%%%%%%%%%%%%%%%%%%%%%%%%%%%%%%%%%%%%%%%%%%%%%%%%%%%%%%%%%%%%%%%%%%%%%

\subsection{General properties for arbitrary $k$}

We may regard $F$ as a continuous piecewise affine map of the $(k-1)$-dimensional simplex
$$
0 \leq x_1 \leq x_2 \leq \dots \leq x_{k-1} \leq 1
$$
into itself. (Although the boundaries 0 and 1 are identified in the original flow, in the
analysis here, we consider them as being distinct points for convenience.)

On the edges of the simplex, $F$ has relatively simple dynamics. Indeed, if
initially all coordinates are equal, then they must all reach the boundary 1
simultaneously. In other words, on the diagonal ($x_i=x$ for all $i$), we have
$F(x,\dots,x)=(t_1,1,\dots,1)$ where $t_1$ depends on $r,s$ and $x$ (for $x=0$,
we have $t_1=1$ independently of $r$ and $s$). Moreover, starting with $x_{k-1}=1$
implies $t_1=0$ which yields
$$
F(x_1,\dots,x_{k-2},1)=(0,x_1,\dots,x_{k-2})
$$
whatever the remaining coordinates $x_1,\dots,x_{k-2}$ are. As a consequence, the edge
$$
\left\{(x,1,\dots,1)\ :\ x\in [0,1]\right\}
$$
is mapped onto
$$
\left\{(0,x,1,\dots,1)\ :\ x\in [0,1]\right\}
$$
which is mapped onto $\left\{(0,0,x,1,\dots,1)\ :\ x\in [0,1]\right\}$ and so on,
until it reaches the edge $(0,\dots,0,x)$, which is mapped back onto the diagonal (after $k$ iterations).

A particular orbit on the edges is the $k$-periodic orbit passing the vertices, and which
corresponds to the single cluster of velocity 1 in the original flow, namely
\[
(0,\dots,0)\mapsto (1,\dots,1)\mapsto (0,1,\dots,1)\mapsto (0,0,1,\dots,1)\mapsto\dots\mapsto (0,\dots,0,1)\mapsto (0,\dots,0).
\]
Geometrically, the corners of the simplex are cyclically permuted by the map $F$. It follows that
the $k-2$ dimensional simplexes (faces) that make up the boundary of the $k-1$ simplex are also cyclically
permuted by $F$. This implies that $F$ cannot have a fixed point on the boundary.
Propositions~\ref{pos_cluster_stability} and \ref{neg_cluster_stability} tell us this orbit 
(which represents a single isolated cluster) must be asymptotically
stable for positive feedback and unstable for negative. Since there are two points of a  periodic
orbit at each boundary of every edge (which are themselves globally $k$-periodic 1-dimensional sets),
and since these points are either both stable or both unstable,
there must be at least one other $k$-periodic orbit on the
edges with coordinates  between 0 and 1. Whether this orbit is unique might depend on parameters.

{\sf Proof of Proposition 2.1.}
Since the simplex is a convex and compact invariant set under $F$, the Brouwer fixed 
point theorem implies the existence of a fixed point. Since the boundary cannot
contain any fixed point, the fixed point is in the interior. Note that the $k$-cluster 
cyclic solutions are fixed points of $F$ and vice versa.
The uniform solution follows by taking $k=n$.
\hfill $\Box$

%%%%%%%%%%%%%%%%%%%%%%%%%%%%%%%%%%%%%%%%%%%%%%%%%%%%%%%%%%%%%%%%%%%%%%%%%%%%%%%%%%
%%%%%%%%%%%%%%%%%%%%%%%%%%%%%%%%%%%%%%%%%%%%%%%%%%%%%%%%%%%%%%%%%%%%%%%%%%%%%%%%%%
%%%%%%%%%%%%%%%%%%%%%%%%%%%%%%%%%%%%%%%%%%%%%%%%%%%%%%%%%%%%%%%%%%%%%%%%%%%%%%%%%%
%%%%%%%%%%%%%%%%%%%%%%%%%%%%%%%%%%%%%%%%%%%%%%%%%%%%%%%%%%%%%%%%%%%%%%%%%%%%%%%%%%

\section{Dynamics for $k=2$}

In this section we study the
dynamics of $2$ cluster systems for the model (\ref{eqn:fmodel}). 
Studying the behavior in the cases
of a small number of clusters is not just a matter of convenience, but
is important from the perspective of applications since presumably only
a small number of clusters can form (for fixed $n$, more clusters implies that
each cluster contains fewer cells and thus can exert less influence) and 
be observable (smaller clusters would produce smaller oscillation in metabolites and
other chemical agents).  In the experiment reported in Figure~\ref{ThePlot} there are $2$ clusters.

\subsection{The map $F$}
 Consider (\ref{eqn:fmodel}) where $f$ is a monotone function.
In the case $k=2$, since only one cluster can exert feedback on the other, RS-feedback (\ref{eqn:fmodel}) simplifies to:
$$
\frac{dx_i}{dt} = \left\{\begin{array}{ccl}
                   1+f(\frac{1}{2}) 
                                   &\text{if}&  x_i  \in [r,1) \text{ and } x_j \in [0,s), j \neq i \\
                   1 &  & \text{otherwise}
              \end{array}\right.  .
$$
Let $\alpha = f(\frac{1}{2})$ for notational simplicity.

For $k=2$, $F$ is defined on the interval $[0,1]$ and is determined by $x_0(t_1)=t_1$  where $t_1$ 
is the time at which $x_1(t)$ reaches $1$.
When regarded as a function of $x_1$ only, its explicit form depends on the parameters $r$ and $s$.
There are two cases depending on the relative sizes of the signaling and responsive regions,
specifically on the size of $(1 +\alpha)s$ with respect to $1-r$.
We have put the details of the computation and analysis of $F(x_1)$ in Appendix~B. 

In the case where $r+(1+\alpha )s<1$, we obtain that $F$ is
a continuous decreasing map:
\[
F(x_1)=\left\{\begin{array}{ccl}
          1-x_1                                            &\text{if}                   &0\leq x_1\leq r-s\\
          1-(1+\alpha )x_1+\alpha (r-s)    &\text{if}   &r-s<x_1< r\\
          1-x_1-\alpha s                          &\text{if}                    &r \leq x_1< 1-(1+\alpha s\\
          \frac{1}{\alpha +1}(1-x_1)                          &\text{if}                    &1-(1+\alpha)s< x_1\leq 1.
\end{array}\right.
\]
In the case where $r+(1 + \alpha )s\geq 1$, the map $F$ is as follows:
\[
F(x_1)=\left\{\begin{array}{ccl}
1-x_1&\text{if}&0\leq x_1< r-s\\
1-(1+\alpha )x_1+\alpha (r-s)            &\text{if}
                        &r-s < x_1\leq \frac{1}{\alpha +1} +\frac{\alpha }{\alpha +1}r - s\\
r-x_1+\frac{1}{\alpha +1}(1-r)        &\text{if}  
                        &  \frac{1}{\alpha +1} +\frac{\alpha }{\alpha +1}r - s   <x_1\leq r\\
\frac{1}{\alpha  +1 } (1-x_1)        & \text{if}  &r< x_1\leq 1.
\end{array}\right.
\]

Calculating the full return map $F^2$ is prohibitively complicated, but in Appendix~B we use these two possible
forms of $F$ to analyze the dynamics. In the next section we summarize the results.

\subsection{Analysis of the dynamics}

\begin{figure}[h!]
\centering
\includegraphics[width=13cm]{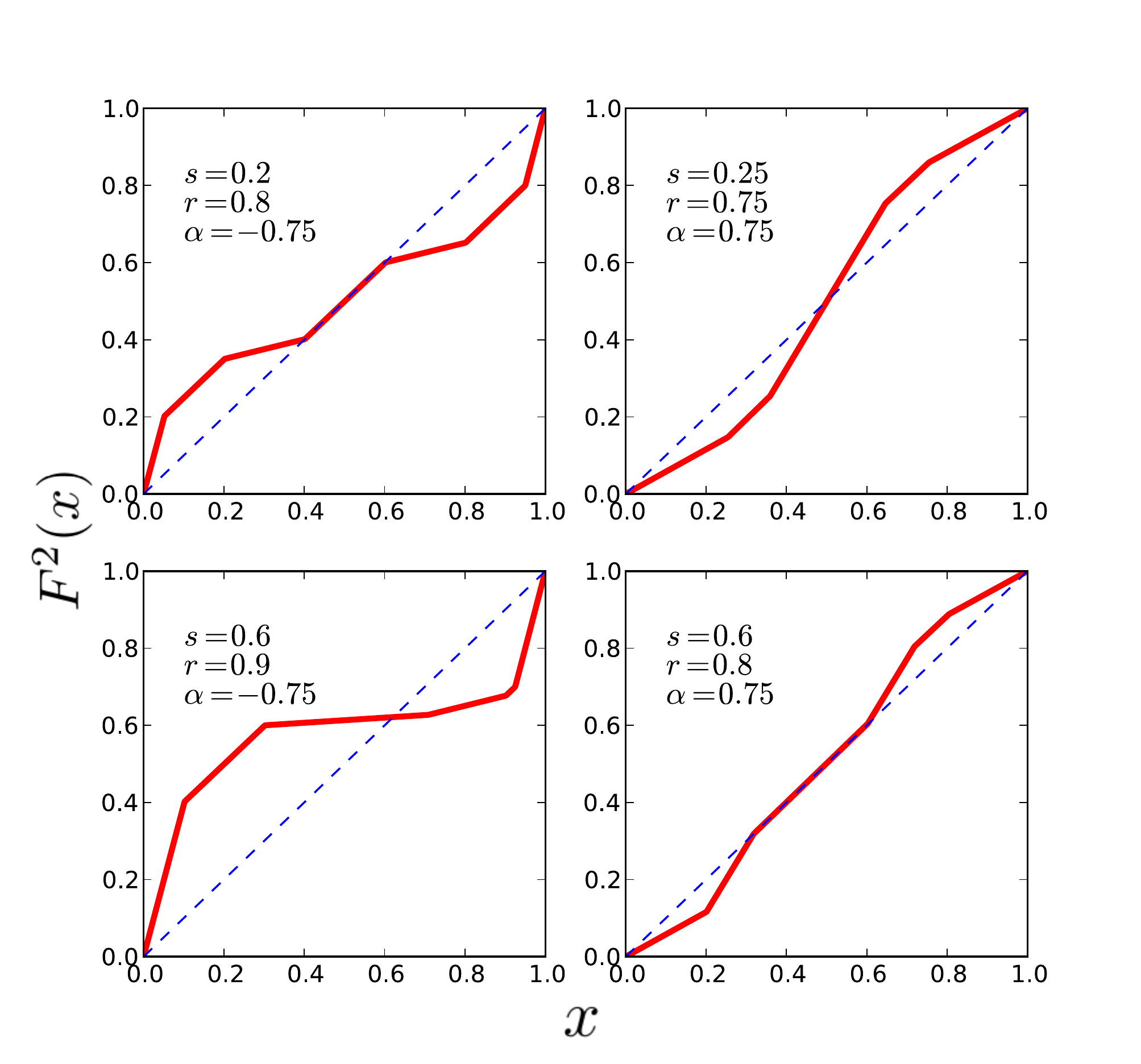}
\caption{Plots of the Poincar\'{e} map $F^2$ in the case $k=2$ for various parameter values ($\alpha = f(\frac{1}{2})$).  Clockwise from top left: interval of fixed points under negative feedback, unstable fixed point under positive feedback, interval of fixed points under positive feedback, stable fixed point under negative feedback.}
\label{F2_return}
\end{figure}

In Appendix~B we find only four distinct types of dynamics; two for positive feedback and two for negative.
\begin{itemize}
\item Positive feedback:
\begin{enumerate}
  \item There is a unique unstable fixed point for $F^2$.
  \item There is an interval of fixed points for $F^2$. 
\end{enumerate}
All other orbits are asymptotic  to the boundary (merger of the two clusters).
\item Negative feedback: 
\begin{enumerate}
  \item There is a unique stable fixed point for $F^2$.
  \item There is an interval of fixed points for $F^2$. 
\end{enumerate}
All other orbits, except the boundary points, are asymptotic to the stable fixed point or the interval of fixed points.
\end{itemize}
These possibilities, for some specific parameter values, are illustrated in Figure~\ref{F2_return}.

An interval of fixed points, we observe in the Appendix B,  can occur not only because the two clusters
may be isolated from each other,  but also in certain other situations.  Namely if either:
\begin{itemize}
\item $x_0$ is in $S$ for the entire time that $x_1$ is in $R$, or,
\item $x_1$ is in $R$ for the entire time that $x_0$ is in $S$,
\end{itemize}
then the unique fixed point of $F$ is neutral and contained in an interval of neutral period 2 points
(fixed points of $F^2$).

In \cite{FY10} we present similar computations for a  subset of parameter values with three clusters 
and positive linear feedback. There we  compute $F$ where $F^3$ is the full return map.
The results there are similar to those reported here; for all the cases examined the three cluster cyclic 
solution is either unstable or in the interior of a set of  neutral periodic solutions (period 3 points of $F$). 
No other periodic orbits were detected and all other initial conditions tend to two cluster or one 
cluster periodic solutions (on the boundary of the domain of the map $F$.

A summary of our studies of $k=2$ and $k=3$ cluster systems is as 
follows. If the system has positive feedback,  then many initial conditions lead to a single cluster,  but
if the initial condition begins with $2$  or $3$ clusters, or close to such, then these clusters 
might persist depending  on  the parameters and initial conditions ({\em e.g.}\ if the clusters are isolated).
If the system has negative feedback, then there may be solutions with $2$ or $3$ clusters (depending on the 
parameter values) that are stable within the set of clustered solutions. One cluster is never stable
under negative feedback. Biologically,  synchrony is likely to appear
in systems with positive feedback  and clustering in systems that have negative
feedback.

%%%%%%%%%%%%%%%%%%%%%%%%%%%%%%%%%%%%%%%%%%%%%%%%%%%%%%%%%%%%%%%%%%%%%%%%%%%%%%%%%%
%%%%%%%%%%%%%%%%%%%%%%%%%%%%%%%%%%%%%%%%%%%%%%%%%%%%%%%%%%%%%%%%%%%%%%%%%%%%%%%%%%
%%%%%%%%%%%%%%%%%%%%%%%%%%%%%%%%%%%%%%%%%%%%%%%%%%%%%%%%%%%%%%%%%%%%%%%%%%%%%%%%%%
%%%%%%%%%%%%%%%%%%%%%%%%%%%%%%%%%%%%%%%%%%%%%%%%%%%%%%%%%%%%%%%%%%%%%%%%%%%%%%%%%%

%
%
%
%\section{Dynamics for $k=3$}
%
%

%%%%%%%%%%%%%%%%%%%%%%%%%%%%%%%%%%%%%%%%%%%%%%%%%%%%%%%%%%%%%%%%%%%%%%%%%%%%%%%%%%
%%%%%%%%%%%%%%%%%%%%%%%%%%%%%%%%%%%%%%%%%%%%%%%%%%%%%%%%%%%%%%%%%%%%%%%%%%%%%%%%%%
%%%%%%%%%%%%%%%%%%%%%%%%%%%%%%%%%%%%%%%%%%%%%%%%%%%%%%%%%%%%%%%%%%%%%%%%%%%%%%%%%%
%%%%%%%%%%%%%%%%%%%%%%%%%%%%%%%%%%%%%%%%%%%%%%%%%%%%%%%%%%%%%%%%%%%%%%%%%%%%%%%%%%

\section{Cyclic $M+1$ Cluster Solutions}

Again consider the model (\ref{eqn:fmodel}) of RS feedback.
Recall that $M = \lfloor (|R|+|S|)^{-1}\rfloor$ is the maximum number of clusters that can exist without mutual interactions.  
In this section we consider the cyclic solutions consisting of $k = M+1$ clusters, with coordinates $ x_0, ..., x_{k-1}$, in the 
dynamics corresponding to (\ref{eqn:fmodel}).  

We again denote the signaling region $S$ by $[0, s)$ and the 
responsive region $R$ by $[r, 1)$.
\begin{prop} \label{prop:Mplus1}
Consider RS feedback of the form (\ref{eqn:fmodel}). 
For any $0 < s < r < 1$, there is a cyclic solution consisting of 
$k = M + 1$ equal clusters of the form $x_0 = 0, \ x_1 = d, \ ..., \ x_{k-1} = (k-1)d$, for some $d>0$.  
Denote $\beta = f(\frac{1}{k})$. If
\begin{equation} 
s < \frac{1}{k}\left(\frac{1 + \beta r}{1 + \beta}\right) \ \ \ \mbox{and} \ \ \ r > \frac{k-1}{k}(1-s\beta   ) \label{sr_ineqs} 
\end{equation}
then the fixed point is unstable for positive $\beta$ and stable  for negative $\beta$. 
 Otherwise, the solution is neutrally stable. Both stability results are within the set of $k$-cluster solutions.
\end{prop}
 
(Note that for $k=2$ the first (resp.\ second) inequality in (\ref{sr_ineqs}) corresponds to the requirement for the second line segment in Figure \ref{RETURN2}b (resp.\ \ref{RETURN2}a) to meet the diagonal, as $k= M+1 =2$ forces $r-s < \frac{1}{2}$)

\begin{figure}[th]
 \centering
 \includegraphics[width = 9cm]{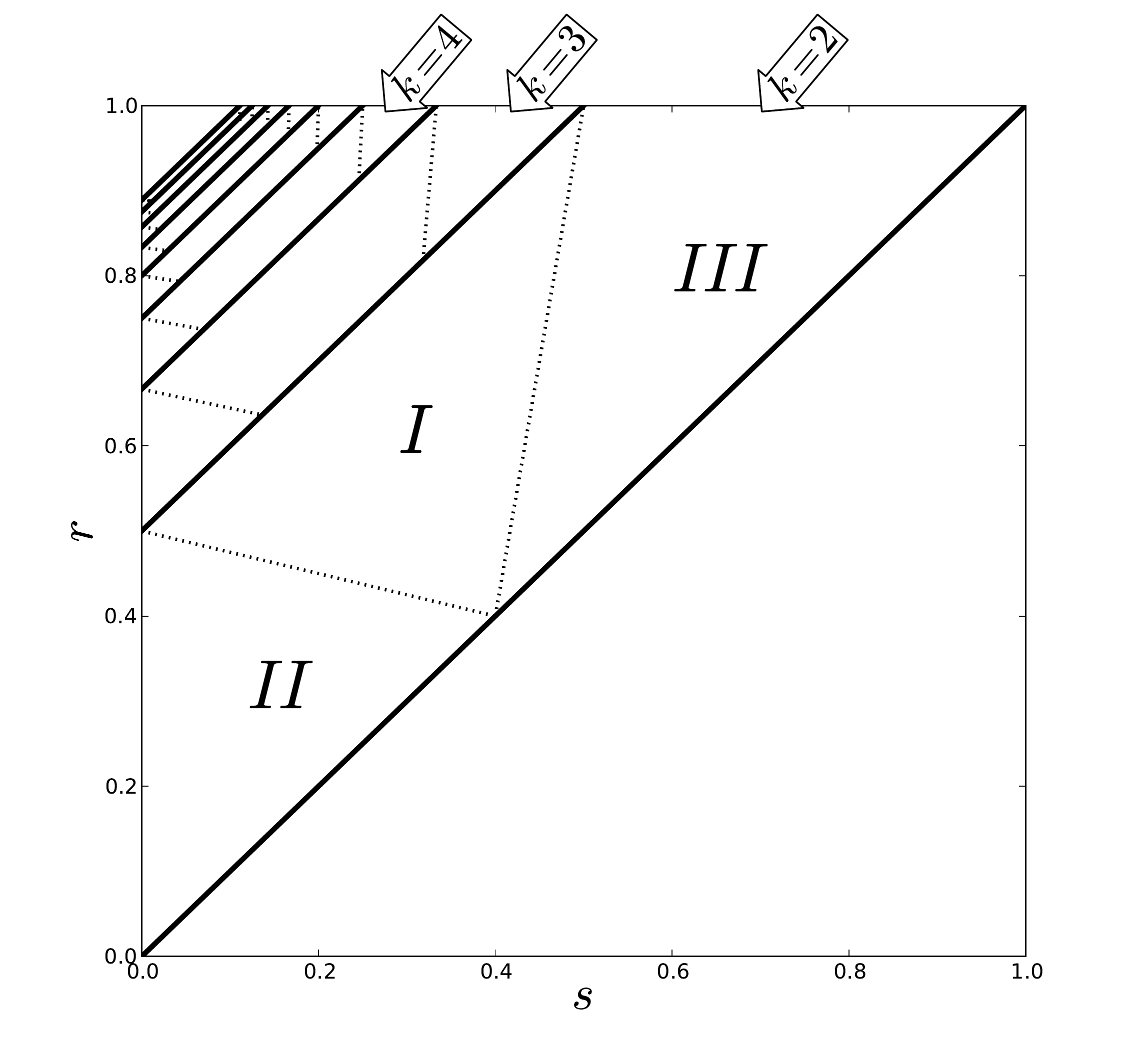}
 \caption{Regions of parameter space for the $k$ cyclic solutions with $k = M+1$. Each diagonal band 
 contains the parameters for a specific $k$ and  is partitioned into three cases.
 In case~I the $k$-cyclic solution is unstable for
 positive feedback and stable for negative.  In cases~II and III the $k$-cyclic solution is neutral and is 
 contained in a set of neutral period $k$ solutions.  
 In this plot  $\beta = f(\frac{1}{k})$ is taken to be $1/k$.}
\label{k-cycles}
\end{figure}

Note that if the solution is neutrally stable, then
the fixed point of $F$ must be contained in the interior of a set of period $k$ points that are also neutrally
stable, since the return map is piece-wise affine. 
Note that this is consistent with the results for two  and three cluster systems.

The proof of Proposition~\ref{prop:Mplus1} is contained in Appendix C.

%%%%%%%%%%%%%%%%%%%%%%%%%%%%%%%%%%%%%%%%%%%%%%%%%%%%%%%%%%%%%%%%%%%%%%%%
%%%%%%%%%%%%%%%%%%%%%%%%%%%%%%%%%%%%%%%%%%%%%%%%%%%%%%%%%%%%%%%%%%%%%%%% 
%%%%%%%%%%%%%%%%%%%%%%%%%%%%%%%%%%%%%%%%%%%%%%%%%%%%%%%%%%%%%%%%%%%%%%%%
%%%%%%%%%%%%%%%%%%%%%%%%%%%%%%%%%%%%%%%%%%%%%%%%%%%%%%%%%%%%%%%%%%%%%%%%
%%%%%%%%%%%%%%%%%%%%%%%%%%%%%%%%%%%%%%%%%%%%%%%%%%%%%%%%%%%%%%%%%%%%%%%%

\section{Discussion}

For both positive and negative feedback, the maximum number of non-interacting clusters, $M$, 
(see Definition~\ref{def:M}), plays a large role in the dynamics of the system.

Our main result for positive feedback  is that the strictly isolated  clustered solutions are neutrally stable
(within the set of $k$ clustered solutions, $k \le M$), 
while the {\em set} of strictly isolated  clustered solutions is locally asymptotically stable (in the full space).
The basin of attraction of this set  extends beyond a small neighborhood of the set.
For (\ref{eqn:fmodel}) with positive feedback and some regions of parameters we have
 proved that the periodic solution of $k=M+1$  equal clusters is completely unstable. 
 In simulations all initial conditions seem to lead to clustered solutions with $k \le M$.
Given that these are only neutrally stable within the set of clustered solutions, solutions may quickly 
approach clustered solutions, but small perturbations of the model may eventually cause the clusters 
to merge into one cluster (full synchronization). 

For negative feedback we observed that isolated  clusters are unstable. This implies that clusters need 
interaction with other clusters to remain coherent. 
When clusters form, the number of clusters $k$ is at least $M+1$.
For the model (\ref{eqn:fmodel}) and under some restrictions on the parameters, we have proved that the 
k-cyclic periodic solution of $k=M+1$ clusters is stable. In numerical experiments, we observed
that clusters sometimes do not form under negative feedback. This perhaps indicates that complicated
bifurcations occur in between regions of clustering.

From simulations in the realistic case $n >> M$ the uniform  solution seems to be unstable
for both positive and negative feedback. Analytic treatment of stability of uniform solutions is
a challenge for future work.

Perhaps the key observation from the mathematical models is that positive feedback 
typically leads to synchrony, while negative feedback systems tend to clustering. 
In practical terms, if a positive feedback mechanism similar to ours were in effect then one would expect to 
encounter synchronization, at least some of the time.
In yeast experiments a single cluster (fully synchronized behavior) seems to be impossible to sustain 
in the laboratory~\cite{breeden,tpb,walker} and in YAO experiments there is clearly no
CDC synchrony. This creates a strong suspicion that the causal mechanism underlying
clustering and YAO involves negative feedback rather than positive.

Acetaldehyde \cite{bier}, ethanol, oxygen and other  membrane-permeating  
metabolites \cite{hjo} have been conjectured as possible mediators of the YAO.
It is natural that these same substances should be considered as signaling agents involved in 
CDC clustering as well.  We have investigated the role of acetaldehyde and ethanol~\cite{PNAS08}, 
and have shown that we can reset phase and induce autonomous oscillation with pulses of acetaldehyde.
Further, we showed that both acetaldehyde and ethanol injections act differentially in different
phases of the CDC and YAO and both act as agents of delay, supporting both our model and 
the role of these substances in processes.

Finally, we note that our analysis is quite general and may find applications outside
the realm of the yeast. We note that other micro-organisms 
can use signaling and coordination of various processes \cite{chen,danino10,dunny,hornby,lyon}. 
It is well recognized that feedback can produce synchrony in physical and biological systems
and this has been extensively studied~\cite{ADGW,Erm,hale98, KKK,mirollo90,PZRO}
(and many others). Clustering is a far less understood phenomenon.
As early as 1977 a model
of the CDC was introduced with periodic blocking at division \cite{rotenberg}, and clustering
behaviour was also observed there.  
 We note that clustering has been observed in Kuramoto type phase oscillators \cite{orosz07} and in 
 certain all-to-all coupled networks of integrate and fire
oscillators \cite{vreeswijk96} and recently a rigorous proof of clustering
and the stability of such solutions was obtained \cite{mauroy08}.
In that work the mechanism that leads to clustering has similarities with the
idealized models in~\cite{BSGY}.

\medskip

{\bf Authors' contributions:}\\
T.Y.~wrote \S 1-4,  worked on \S 2-6 and coordinated the work and manuscript.
B.F.~worked on \S 2, 3, 5, 6, wrote \S 5, 6 and produced graphics.
R.B.~worked on and wrote \S 7 and produced graphics and simulation code.
G.M.~contributed to the final form of \S 6 and \S 7.
E.B.~conceived the project, worked on sections 4, 7 and edited the manuscript.

{\bf Acknowledgments:}\\
B.F.~thanks the Courant Institute (NYU) for hospitality. He~was supported by CNRS
and by the EU Marie Curie fellowship PIOF-GA-2009-235741.
E.B., T.Y.\ and this work were supported by the NIH-NIGMS grant R01GM090207.
The authors thank the referees for invaluable corrections and comments 
that greatly improved the manuscript.

\appendix

\section{Appendix - Relation to PDE models}

Consider the hyperbolic conservation law:
\begin{equation}\label{pde}
      \frac{\partial u}{\partial t} +  \frac{\partial}{\partial x} \left( b(x,[u]) u \right) = 0
\end{equation}
for $x$ on the unit circle, where $u(x,t) \ge 0$ represents the density of cells of size $x$
and $[u]$ indicates that $b$ depends functionally on the solution $u$ \emph{i.e.} $b(x, [u])$ 
depends on the values of $u$ at all points $x \in S^1$ at time $t$. This is a conservation law
since $\int_{S^1} u(s,t) \, dx$ is a constant in time.
This is the natural PDE version of our model of RS feedback  
(compare \cite{diekmann84,hannsgen85,henson,rotenberg}) and we provide it here for comparison
with the existing literature. 

If $b$ has the form of RS feedback, then given any
constant $c>0$, there is a solution $u(x)$ such that $u(x) \equiv c$ for $x \in R^c$ and
$u(x)$ for $x \in R$ is determined explicitly by  the feedback $b(x,[c])$, {\em i.e.},
\begin{equation}\label{steadyu}
  u(x) = \begin{cases}
            c, \quad \textrm{for} \quad x \notin R,\\
          \frac{c}{b(x,[c])}  , \quad \textrm{for} \quad x \in R.
         \end{cases}
\end{equation}
This solution is analogous to the uniform solution in Proposition~2.3.

Note that the corresponding version of (\ref{eqn:fmodel}) is given by:
\begin{equation}\label{pde_linear}
    b(x,[u]) = \begin{cases}
                      1, \quad \textrm{if} \quad x \notin R \\
                      1 + f(I), \quad \textrm{if} \quad x \in R,
                     \end{cases}
\end{equation}
where
$$
    I =\int_S u \, dx,
$$
{\em i.e.}\ the fraction of cells in $S$.
For this feedback term, the uniform solution is piecewise constant with value $c/(1 + f( c |S|))$ in $R$.
    
We are not aware of any treatment of PDEs such as (\ref{pde}) with $b$ of the form $b(x,[u])$.
When there is no CDC feedback and diffusion is added to the PDE model (in various ways), then the uniform solution
is asymptotically stable \cite{diekmann84,hannsgen85}.

%%%%%%%%%%%%%%%%%%%%%%%%%%%%%%%%%%%%%%%%%%%%%%%%%%%%%%%%%%%%%%%%%%%%%
%%%%%%%%%%%%%%%%%%%%%%%%%%%%%%%%%%%%%%%%%%%%%%%%%%%%%%%%%%%%%%%%%%%%%
%%%%%%%%%%%%%%%%%%%%%%%%%%%%%%%%%%%%%%%%%%%%%%%%%%%%%%%%%%%%%%%%%%%%%
%%%%%%%%%%%%%%%%%%%%%%%%%%%%%%%%%%%%%%%%%%%%%%%%%%%%%%%%%%%%%%%%%%%%%

\section{Appendix - Analysis for $k=2$}

Recall the notation $\alpha = f(\frac{1}{2})$.

We report the details of computations in the case where $(1 + \alpha )s<1-r$. 
 The other case can
be treated similarly and we only give below the resulting expression of $F$. When 
$(1 + \alpha )s<1-r$,
there are 4 situations depending on the location of $x_1$:
\begin{itemize}
\item[$\bullet$] $x_1\leq r-s$. In this case, $x_0$ leaves $S$ before $x_1$ enters $R$.
The point $x_1$ is not submitted to any feedback and hence $x_1(t)=x_1+t$ for all $t$
which implies $F(x_1)=t_1=1-x_1$. (The occurrence of this case is independent 
of $r+(1 + \alpha )s <1$.)
\item[$\bullet$] $r-s<x_1\leq r$. Here $x_1$ is influenced, but not during the
entire responsive region since $x_0$ gets out of $S$ before $x_1$ reaches 1. More precisely, we
have\footnote{the occurrence of the late phase where $x_1$, although being in $R$, moves
with velocity 1, occurs due to the condition $r+(1 + \alpha )s<1$.}
\[
x_1(t)=\left\{\begin{array}{ccl}
x_1 + t                     &\text{if}&0<t<r-x_1\\
r+(1+\alpha )(t-r+x_1)&\text{if}&r-x_1<t<s\\
r+(1+\alpha )(s-r+x_1)+t-s&\text{if}&s\leq t
\end{array}\right.
\]
It follows that $F(x_1)=1-(1+\alpha )x_1+\alpha (r-s)$.
\item[$\bullet$] $r<x_1\leq 1-(1+\alpha )s$. Then 
\[
x_1(t)=\left\{\begin{array}{ccl}
x_1+(1+\alpha )t          &\text{if}&0<t<s\\
x_1+(1+\alpha )s+t-s                 &\text{if}&s\leq t
\end{array}\right.
\]
from which we obtain $F(x_1)=1-x_1-\alpha s$.
\item[$\bullet$] $1-(1+\alpha )s<x_1\leq 1$. In this case, $x_1$ starts sufficiently close
to 1 to have velocity $1+\alpha $ when reaching the boundary. We have 
$F(x_1)=\frac{1}{\alpha +1}(1-x_1)$.
\end{itemize}

\begin{figure}[t!]
\centering
\includegraphics[width=6cm,height=8cm]{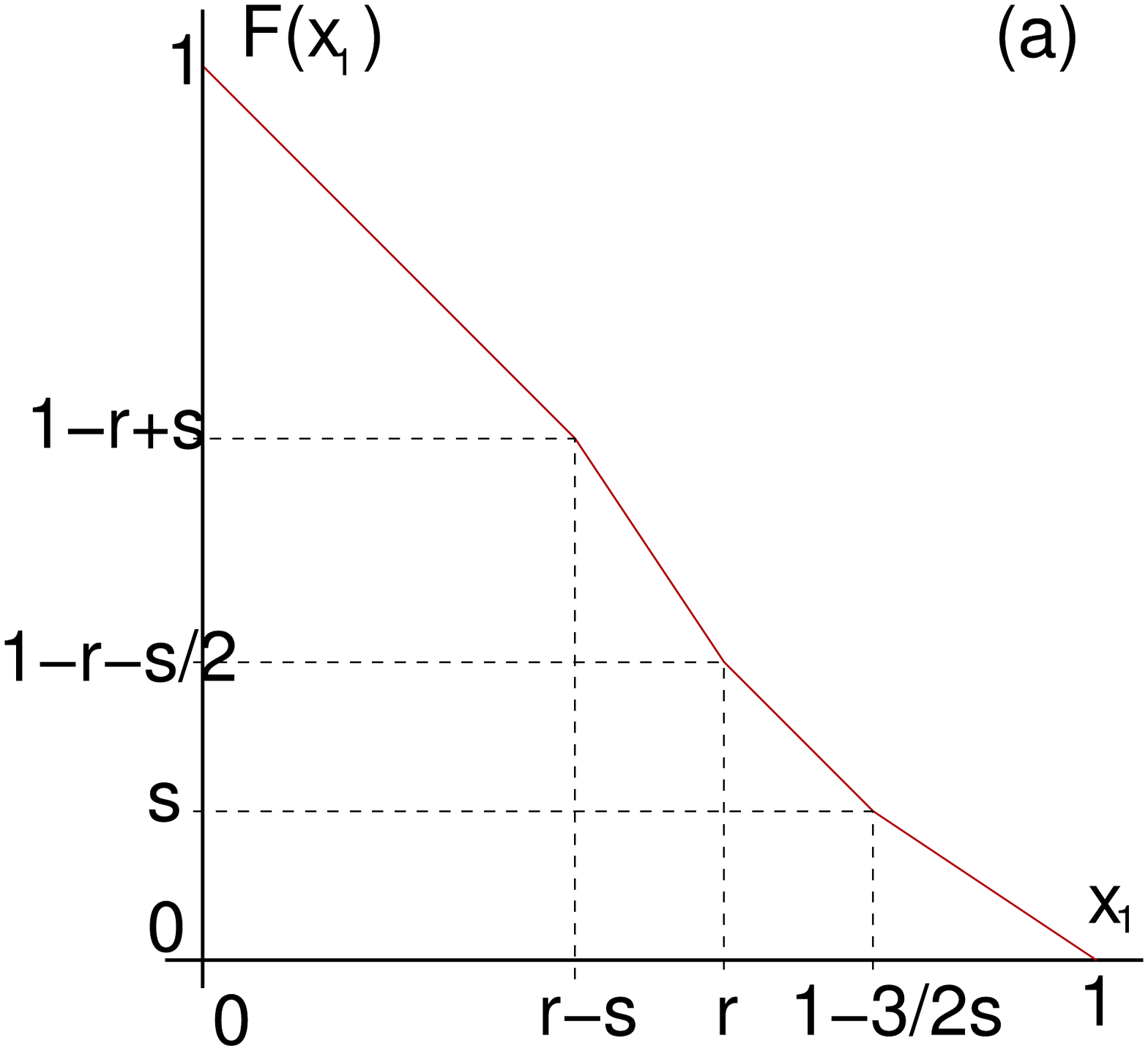}
\hspace*{.5cm}
\includegraphics[width=6cm,height=8cm]{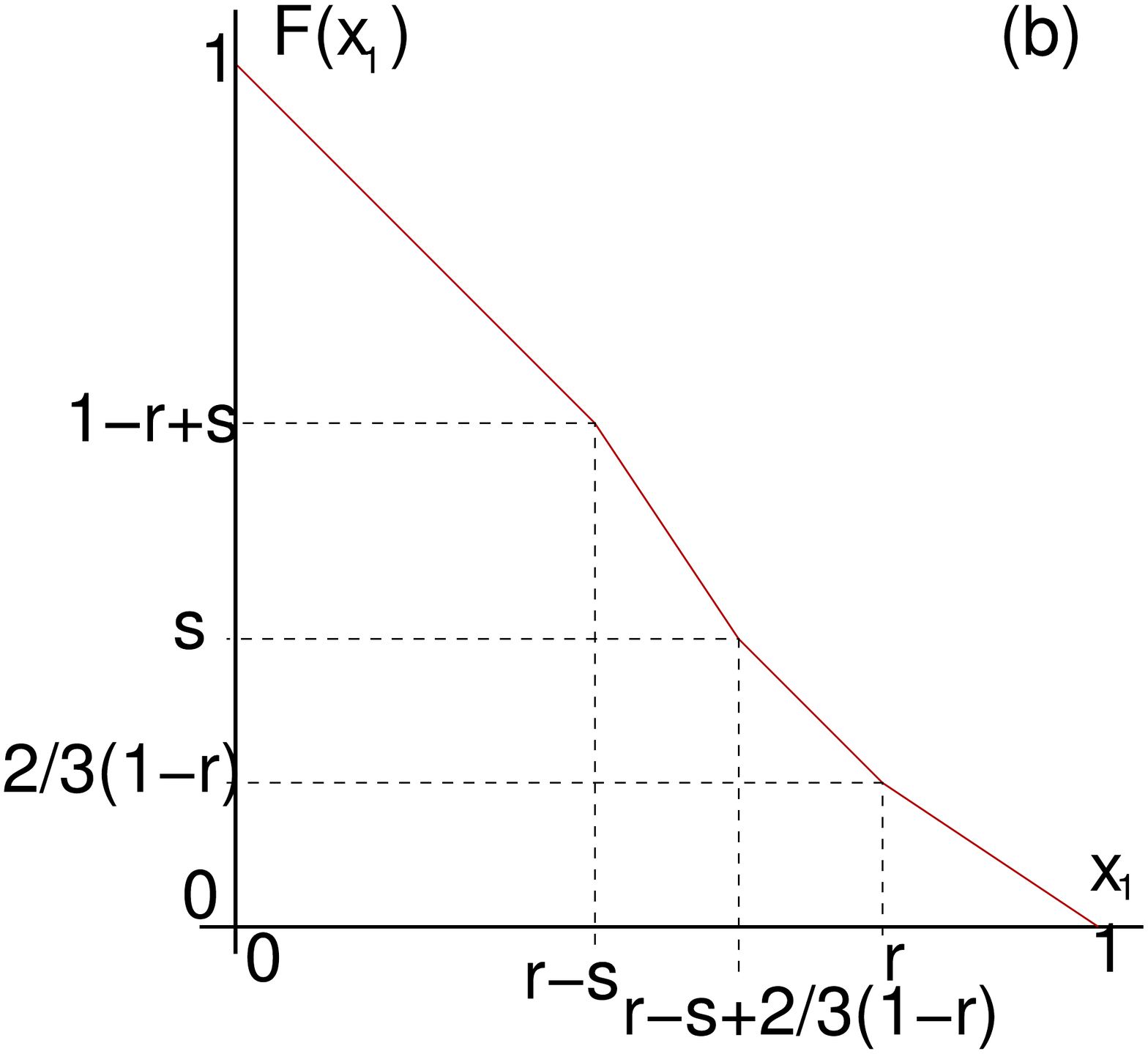}
\caption{Plots of the mapping $F$  for $k=2$ and $\alpha = f(\frac{1}{2}) = \frac{1}{2}$. (a) $r+\frac{3}{2}s<1$. (b) $r+\frac{3}{2}s\geq 1$.
If $r-s >1/2$ (which implies $M\ge 2$) then there is a neutral fixed point at $1/2$ that represents an isolated $2$ 
cluster cyclic solution. If the second segment intersects the diagonal, then there is an isolated $2$ cluster cyclic solutions
that is stable for $\alpha  <0$ and unstable for $\alpha  >0$. There also may exist neutral fixed points for $r-s<1/2$ ($M=1$) and certain conditions on 
the parameters where the third piecewise segment of $F$ intersects the diagonal line. These fixed points represent
$2$ cluster cyclic solutions that are not isolated, but yet are neutral.}
\label{RETURN2}
\end{figure}

Recapitulating, in the case where $r+(1+\alpha )s<1$, we obtain the following expression
of a continuous decreasing map whose plot is given in Figure \ref{RETURN2} (a)
\begin{equation}
F(x_1)=\left\{\begin{array}{ccl}
          1-x_1                                            &\text{if}                   &0\leq x_1\leq r-s\\
          1-(1+\alpha )x_1+\alpha (r-s)    &\text{if}   &r-s<x_1< r\\
          1-x_1-\alpha s                          &\text{if}                    &r \leq x_1< 1-(1+\alpha) s\\
          \frac{1}{\alpha +1}(1-x_1)                          &\text{if}                    &1-(1+\alpha )s< x_1\leq 1.
\end{array}\right.
\label{eqn:Fcase1}
\end{equation}
In the case where $r+(1 + \alpha )s\geq 1$, the map $F$ is as follows (see Figure \ref{RETURN2} (b)).
\begin{equation}
F(x_1)=\left\{\begin{array}{ccl}
1-x_1&\text{if}&0\leq x_1< r-s\\
1-(1+\alpha )x_1+\alpha (r-s)            &\text{if}
                        &r-s < x_1\leq \frac{1}{\alpha +1} +\frac{\alpha }{\alpha +1}r - s\\
r-x_1+\frac{1}{\alpha +1}(1-r)        &\text{if}  
                        &  \frac{1}{\alpha +1} +\frac{\alpha }{\alpha +1}r - s   <x_1\leq r\\
\frac{1}{\alpha  +1 } (1-x_1)        & \text{if}  &r< x_1\leq 1.
\end{array}\right.
\label{eqn:Fcase2}
\end{equation}

As argued for arbitrary $k$, the map $F$ has a  $k$-periodic orbit which, for $k=2$, is
composed of the boundaries 0 and 1. By the Intermediate Value Theorem, it must also have a
fixed point on the diagonal.

The graph of $F$ coincides with the anti-diagonal ($1-x_1$) for $x_1\leq r-s$.
For $\alpha = f(\frac{1}{2}) >0$ it is strictly lower than this line if $r-s<x_1<1$; but if $\alpha  <0$, then
 it is strictly greater than $1-x_1$ for $r-s<x_1<1$.
The dynamics can be characterized completely for arbitrary parameters when $k=2$. The following conclusions hold:
\begin{itemize}
		\item If $r-s > \frac{1}{2}$, each point of the interval $[1-r+s, r-s]$ is part of a 2-periodic
		           orbit $x \mapsto 1-x \mapsto x$, except for the point $x = \frac{1}{2}$ which is fixed.  
		           The return map $F^2$ thus has an interval of neutrally stable fixed points centered 
		           around $x_1 = \frac{1}{2}$. 
		           \begin{itemize}
		               \item[$\bullet$] If $\alpha >0$, if $x_1$ is above (resp. below) this neutral interval, we 
				    have $F^2(x_1) > x_1$ (resp. $F^2(x_1) < x_1$) and so any initial point converges to 1 (resp. 0).
		               \item[$\bullet$] If $\alpha <0$, if $x_1$ is above (resp. below) this neutral interval, we have 
		                        $F^2(x_1) < x_1$ (resp. $F^2(x_1) > x_1$) and so any initial point converges to the boundary of 
		                        the neutral interval. 
	                     \end{itemize}           	
	          	\item If $r-s = \frac{1}{2}$, there is a unique fixed point. It is  stable  for negative $\alpha $ 
		                                       and unstable for positive $\alpha $. 
		\item If $r-s<\frac{1}{2}$ there are three possibilities depending on where the diagonal line $x=y$ intersects
		           the graph of $F$  (see Figure~\ref{RETURN2}).
		          \begin{itemize}
		               \item[$\bullet$] If the diagonal intersects the second segment 
		                           then there is a unique fixed point which is stable if $\alpha <0$ 
		                           and unstable if $\alpha$ is positive.  
		               \item[$\bullet$] If the diagonal intersects the third segment of $F$, then there is again an interval
		                                          of neutral period 2 points. The  edge of the interval is stable 
		                                           for negative $\alpha $ and unstable  for   positive $\alpha $. 
		                \item[$\bullet$] If the diagonal hits the boundary between segments 2 and 3 of $F$ then 
		                                          there is a unique fixed point which stable for negative $\alpha $ and 
		                                           unstable for positive $\alpha $. 
		          \end{itemize}       
\end{itemize}
All of these possibilities are summarized in just four distinct types of behaviour in section 6.2 and Figure~\ref{F2_return}.

The condition $r-s > 1/2$ corresponds to $M \ge 2$. In such a case the cyclic
$2$ cluster solution consists of isolated clusters and it is contained in an interval of neutral period
two points. This interval is an attractor for negative feedback and a repeller for positive feedback.

Note that the condition $r-s < 1/2$ corresponds to $|R|+|S| > 1/2$ which implies that
$M = \lfloor (|R|+|S|)^{-1} \rfloor = 1$. Thus when the 2 clusters cannot be isolated, there is a cyclic 2 cluster solution 
which is a fixed point of $F$. This solution may be unique and stable (negative feedback), unique and unstable 
(positive feedback), or neutral, depending on the parameters.

One can easily calculate that the diagonal line cannot intersect the fourth segment of the graph of $F$ in
either  case (\ref{eqn:Fcase1}) or (\ref{eqn:Fcase2}).

In (\ref{eqn:Fcase1}) it is seen that the third segment can be intersected by the diagonal by making $r$
sufficiently small. This corresponds to the $x_1$ coordinate of the fixed point being greater than $r$ so
that  the second cluster begins in the responsive region $R$. From the conditions, it is still in 
$R$ when $x_0$ leaves $S$. There is an interval of fixed points for $F^2$ even though
the clusters are not isolated.

The diagonal also can intersect the third segment for (\ref{eqn:Fcase2}) if $x_0$ is
in $S$ when $x_1$ enters $R$,  and when $x_1$ reaches $1$, $x_0$ is still in $S$.  
Thus we have another case of interacting clusters that still leads to a neutral fixed point.

%%%%%%%%%%%%%%%%%%%%%%%%%%%%%%%%%%%%%%%%%%%%%%%%%%%%%%%%%%%%%%

\section{Appendix  - Proof of Proposition~7.1}

\noindent
{\sc Proof.}
The evolution of the system may be described qualitatively in terms of `milestones,' \emph{e.g.} $x_0$ reaching s, 
or $x_{k-1}$ reaching 1.  We first consider a system which evolves via the following sequence of events (which we call {\bf Case 1}):
\begin{equation} 
     x_{k-1} \mapsto r,      \qquad    x_0 \mapsto s,     \qquad  x_{k-1} \mapsto 1  \label{case_1_sequence} 
\end{equation}
at which point after relabeling $x_i := x_{i-1}$ we have the initial condition again.  By calculating the time taken in each step 
and finding the final value of $x_0$, one can use the relation $x_{k-1} =(k-1)d$ to solve analytically for $d$, getting:
\begin{equation}
       d = \frac{1+\beta(r-s)}{k+\beta(k-1)} \ \mbox{.} \label{case_1_d} 
\end{equation}
Note that the Case 1 sequence will occur provided that $s < x_1 = d$ and $x_{k-1} = (k-1)d < r$.  Using (\ref{case_1_d}) in 
these relations gives (\ref{sr_ineqs}).

Two other sequences of events are possible and occur when we drop each of the constraints in (\ref{sr_ineqs}). {\bf Case 2}, when we allow $ r < (k-1)d$, is characterized by:
\begin{equation} 
      x_0 \mapsto s, \qquad x_{k-1} \mapsto 1, \nonumber 
\end{equation}
and {\bf Case 3}, when we allow $ s > d $, by:
\begin{equation} 
     x_1 \mapsto s,    \qquad  x_{k-1} \mapsto r, \qquad x_{k-1} \mapsto 1  \mbox{.} \nonumber 
     \end{equation}
     
Following the same procedure as in Case 1, we obtain for Case 2:
\[ 
     d = \frac{1-s \beta}{k} 
\]
and for Case 3:
\[
 d = \frac{1 + r \beta}{k(1 + \beta)} \ \mbox{.} 
 \]

A simple calculation shows that for a given $M$, cases 1, 2 and 3 exhaust the parameter set in $(r,s)$.  We observe that there is never more than one cluster in the signaling region when the response region is nonempty, and thus the dynamics of the system are determined entirely by $\beta = f(\frac{1}{k})$.
A graphical representation of the regions of parameter space corresponding to these three cases can be seen in figure (\ref{k-cycles}).

The map $F$ is affine in a neighborhood of the fixed point, i.e.\ $\vec{x} \mapsto A\vec{x} + \vec{b}$ 
where $\vec{x} = (x_0,...,x_{k-1})^T$ and $A$ is a matrix. We next analyze $A$ in the three cases.

{\bf Case 1:}
\[ A = \left[
	\begin{array}{cccccc}
	0 & 0 & 0 & \cdots & 0 & -(1 + \beta) \\
	1 & 0 & 0 & \cdots & 0 & -(1 + \beta) \\
	0 & 1 & 0 & \cdots & 0 & -(1 + \beta) \\
	& \vdots & & \ddots & \vdots & \vdots \\
	0 & 0 & 0 & \cdots & 0 & -(1 + \beta) \\
	0 & 0 & 0 & \cdots & 1 & -(1 + \beta) \\
	\end{array}
\right] \]
Thus we may determine the stability of the fixed point by studying the eigenvalues $A$, which
has characteristic equation 
\begin{equation} 
          -\lambda^{k-1} - (1 + \beta)(\lambda^{k-2} + \lambda^{k-3} + \cdots + \lambda + 1) = 0 \mbox{.}
\label{charEquation}
\end{equation}
Notice that $\lambda =1$ can easily be ruled out as a root.
For \( \lambda \neq 1 \) we can rewrite (\ref{charEquation}) as
\[
\frac{1}{1 + \beta}\lambda^n + \sum_{i = 0}^{n - 1} \lambda^i \ = \
\frac{\lambda^n}{1 + \beta} + \frac{\lambda^n - 1}{\lambda - 1} \ = \ 0 \mbox{.}
\]
After simplification, we see that \( \lambda \neq 1 \) is a solution of (\ref{charEquation}) if and only if
\begin{equation}
\left( \frac{\lambda + \beta}{1 + \beta} \right)\lambda^n = 1 \mbox{.}
\label{rootRequirement}
\end{equation}

Now, suppose \( \beta  > 0 \).  If \( | \lambda | \le 1 \) (and $\lambda \neq 1$), then 
\( | \lambda^n | \le 1 \) and \( | \lambda + \beta | < | 1 + \beta |  \).  Thus
\[
\left| \frac{\lambda + \beta}{1 + \beta}\right| \left| \lambda^n \right| < 1  \mbox{,}
\]
\emph{i.e.} \( \lambda \) cannot satisfy~(\ref{rootRequirement}). Thus for positive feedback all of the eigenvalues 
lie outside the unit disk and so the map is unstable at the fixed point. Further, it is not only unstable, but is 
unstable w.r.t.\ all possible perturbation directions.

For the case \( \beta < 0 \), suppose \( | \lambda | > 1 \).  Write \( \lambda + \beta \) 
as \( \lambda - (-\beta) \).
Then by the reverse triangle inequality,
\( 
     |\lambda + \beta| = | \lambda - (-\beta) | 
           \geq \left| |\lambda | - |-\beta| \right| 
                = \left| |\lambda | + \beta \right| = | \lambda | + \beta > 1 + \beta \).
Thus we have
\[
\left| \frac{\lambda + \beta}{1 + \beta}\right| \left| \lambda^n \right| > 1 \mbox{,}
\]
and~(\ref{rootRequirement}) is not satisfied. Also, if $|\lambda|=1$ but $\lambda \neq 1$, then
$|\lambda +\beta| < 1+\beta$. Therefore if \( \beta < 0 \) then all the eigenvalues 
of \( A \) lie on the interior of the unit disc, and the map is stable. \hfill $\Box$

\begin{figure}[!h]
 \centering
 \includegraphics[width = 15cm]{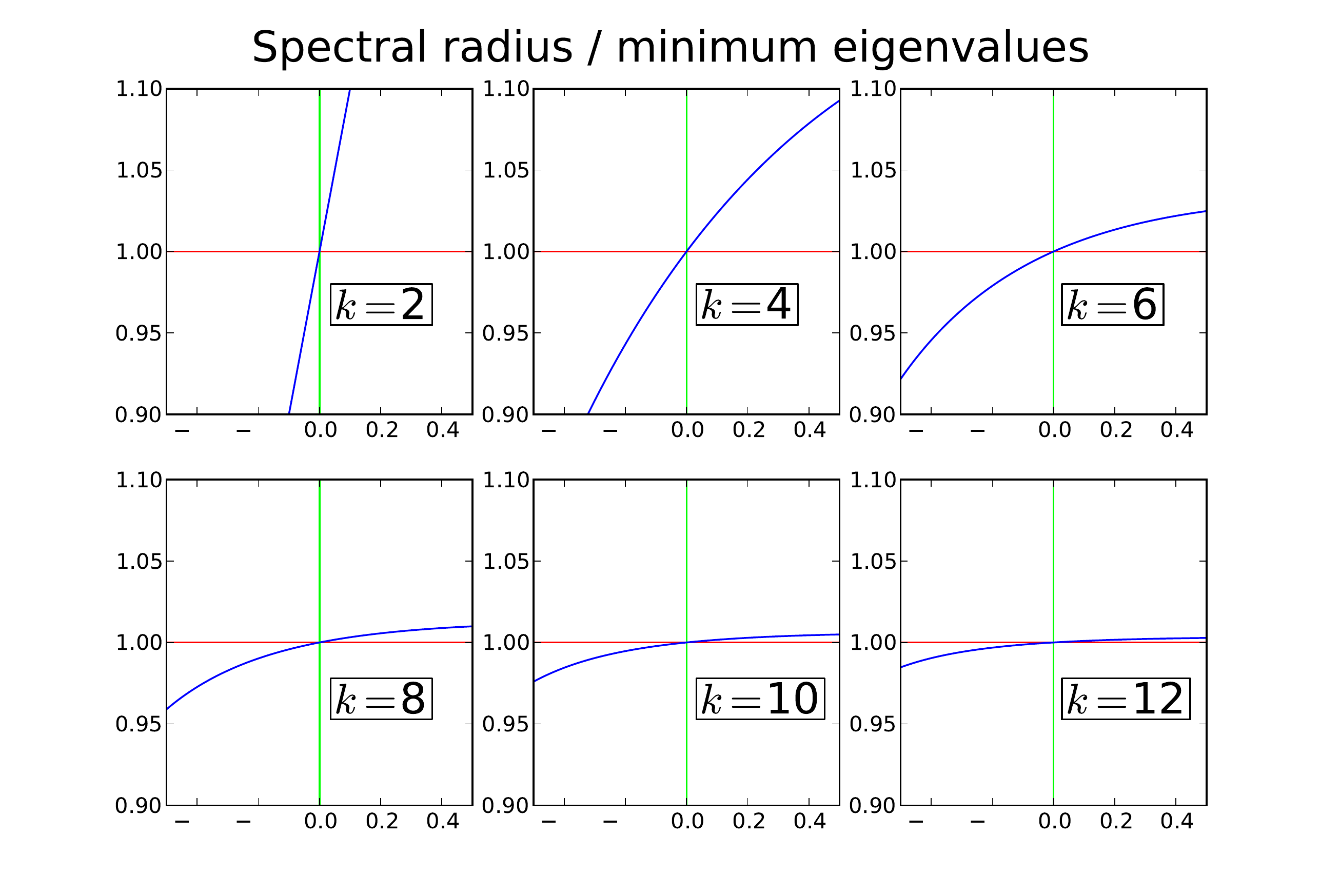}
 \caption{Case 1. The spectral radius / smallest eigenvalue modulus for \( k = 2, 4, \ldots, 12 \).  The \( x \)-axis 
 is the feedback parameter \( \beta \).  For negative \( \beta \), the \( y \)-axis is the spectral radius of \( A \).  
 For positive \( \beta \), the smallest eigenvalue (w.r.t.\ modulus) is plotted. The plots show that the $k = M+1$ 
 solution is stable for negative feedback and unstable (in all directions) for positive feedback. Notice that as 
 $k$ grows, the stability/instability becomes weaker.}
\label{eigs}
\end{figure}

The stability results for case 1 are  illustrated in Figure~\ref{eigs} for \( k = 2, 4, \ldots, 12 \) and \(\beta =  f(\frac{1}{k}) \in (-.5, .5) \).

{\bf Case 2 and Case 3:}\\
In cases 2 and 3, the linear part of the map at the fixed point is represented by the matrix:
\[ 
A = \left[
	\begin{array}{cccccc}
	0 & 0 & 0 & \cdots & 0 & -1 \\
	1 & 0 & 0 & \cdots & 0 & -1 \\
	0 & 1 & 0 & \cdots & 0 & -1 \\
	& \vdots & & \ddots & \vdots & \vdots \\
	0 & 0 & 0 & \cdots & 0 & -1 \\
	0 & 0 & 0 & \cdots & 1 & -1 \\
	\end{array}
	\right] .
\]
This matrix has characteristic equation \( \lambda^{k-1} + \lambda^{k-2} + ... + \lambda + 1 = 0 \) whose roots all 
have absolute value 1.  Thus the map is neutrally stable in both cases 2 and 3.
\hfill $\Box$

\end{document}